\documentclass[final,1p,times]{elsarticle}
\usepackage{ulem}

\usepackage{mdframed}
\usepackage{epsfig}

\usepackage{amsmath,amsthm,verbatim,amssymb,amsfonts,amscd, graphicx}
\usepackage{latexsym,graphicx, epsfig, bm, epstopdf, float, color}
\usepackage{tikz,siunitx}
\usetikzlibrary{arrows.meta,calc,decorations.markings}
\usepackage{scalerel}
\usetikzlibrary{shapes,arrows,patterns}
\usetikzlibrary{shapes.geometric}
\usepackage{pgfplots}
\usetikzlibrary{shapes.geometric,shapes.symbols}
\usetikzlibrary{patterns}
\usetikzlibrary{tikzmark}

\usepackage{caption}
\usepackage{subcaption}
\usepackage{graphics,mathtools,hyperref,cleveref}
\usepackage{dsfont}
\usepackage{framed}
\usepackage[most]{tcolorbox}
\usepackage{xcolor}
\usepackage{multirow}

\newcommand{\av}[1]{\left\{\!\left\{#1\right\}\!\right\}} 
\newcommand{\jump}[1]{\lbrack\!\lbrack #1 \rbrack\!\rbrack} 

\newcommand{\DG}{\text{DG}}

\newcommand{\R}{\mathbb{R}}

\newcommand{\x}{\textbf{x}}

\newcommand{\Tau}{\mathcal{T}}
\newcommand{\vertiii}[1]{{\left\vert\kern-0.25ex\left\vert\kern-0.25ex\left\vert #1 
    \right\vert\kern-0.25ex\right\vert\kern-0.25ex\right\vert}}
\newcommand{\eff}{\mathrm{eff}}
\theoremstyle{remark}
\newtheorem{remark}{Remark}

\usepackage[textsize=footnotesize]{todonotes}

\begin{document}

\begin{frontmatter}



\title{Heat Transfer Modeling in Enhanced Geothermal Energy:\\ A Three-Temperature Approach for Solid, Injected, and\\ Residing Fluids}


\author[FSU]{Yi-Yung Yang} 
\author[FSU]{Sanghyun Lee} 
\author[DK]{Dmitri Kuzmin} 


\affiliation[FSU]{organization={Department of Mathematics, Florida State University},
            addressline={1017 Academic Way}, 
            city={Tallahassee},
            postcode={32306-4510}, 
            state={FL},
            country={USA}}
\affiliation[DK]{organization={Institute of Applied Mathematics, TU Dortmund University},
            addressline={Vogelpothsweg 87}, 
            city={Dortmund},
            postcode={D-44227}, 
            country={Germany}}

\begin{abstract}
Enhanced geothermal systems (EGS) are governed by strongly coupled, advection-dominated flow and heat transfer in fractured porous media. Conventional models typically assume local thermal equilibrium with a single effective fluid temperature or, at best, an averaged pore-fluid temperature, so that the thermal evolution of the injected cold fluid is only inferred indirectly. In this work, we develop a new local thermal non-equilibrium (LTNE) formulation that explicitly resolves the temperature of injected fluid as it travels through the reservoir and exchanges heat with the hot rock and resident fluid. The key modeling ingredient is a concentration variable that tracks the injected fluid and drives a three-way LTNE coupling: rock temperature, resident-fluid temperature, and injected-fluid temperature. This allows us to distinguish, at the continuum level, how newly injected fluid parcels are heated by conductive and convective exchange, and to predict injected-fluid temperature at production wells without relying on bulk averages. To discretize the resulting nonlinear, advection-dominated system, we employ an enriched Galerkin (EG) finite element method for Darcy flow, temperature, and concentration, ensuring local mass conservation with a moderate number of degrees of freedom. A tailored flux-corrected transport (FCT) strategy is constructed for the EG discretization of the concentration and temperature equations to enforce a discrete maximum principle and suppress nonphysical oscillations while preserving local conservation. Time integration is carried out with an IMPES-type splitting combined with a strong-stability-preserving third-order Runge--Kutta (SSP RK2) scheme. Numerical experiments for fractured EGS configurations demonstrate that the proposed LTNE--EG--FCT framework can resolve injected-fluid heating paths and thermal breakthrough behavior that are not captured by standard single-temperature or averaged LTNE models.
\end{abstract}







\end{frontmatter}

\section{Introduction} \label{sec:Introduction}

Enhanced geothermal systems (EGS) \cite{OLASOLO2016133,NATH2024213370} extract heat from deep, low-permeability rock formations by injecting a relatively cold working fluid, circulating it through a fractured porous medium, and producing it at elevated temperatures. 
The efficiency and long-term viability of EGS operations are governed by strongly coupled fluid flow and heat transfer processes in highly heterogeneous media, frequently operating in an advection-dominated regime~\cite{mcclure2014investigation}. 
Accurately resolving how injected cold fluid exchanges heat with the hot rock matrix and resident pore fluid is essential for predicting thermal breakthrough, reservoir lifetime, and power output.

Most geothermal reservoir simulators rely on the local thermal equilibrium (LTE) assumption \cite{QUINTARD19952779}, in which a single effective temperature represents the combined rock–fluid system. 
Even when local thermal nonequilibrium (LTNE) \cite{NIELD19993245,PARHIZI2021120538} models are employed, they typically introduce separate temperatures for the solid and fluid phases while still treating the pore fluid as a single averaged component. 
Consequently, the thermal evolution of injected cold fluid is not explicitly resolved; instead, predictions rely on averaged temperatures within the porous medium. 
{This averaging obscures the detailed heating dynamics of newly injected fluid parcels and limits the ability to resolve multiscale phenomena such as fracture~\cite{vik2018heat,caulk2016experimental,lee2025phase} channeling, fingering instabilities, and localized heat-exchange hot spots.}

{To address this limitation, we develop a three-way local thermal nonequilibrium (three-way LTNE) framework that distinguishes among injected fluid, resident fluid, and solid matrix temperatures.} 
The governing system couples Darcy flow for a single-phase fluid with separate energy equations for the injected fluid, resident fluid, and rock matrix, and incorporates a concentration equation to track injected-fluid transport. 
{This formulation enables explicit tracking of injected-fluid heating along flow paths, thereby resolving thermal processes that are only indirectly represented in conventional LTE or two-temperature LTNE models.}

The resulting system is strongly coupled and advection dominated, posing significant challenges for numerical discretization. 
To obtain physically consistent and locally conservative solutions, we employ an enriched Galerkin (EG) finite element method for both the Darcy flow and transport equations \cite{Kuzmin-EG,LEE2016,LEE201719,yi2024physics,lee2025thermo}. 
The EG formulation ensures elementwise mass conservation while maintaining a relatively low number of degrees of freedom, making it suitable for large-scale geothermal simulations. 
For temporal coupling, we adopt a modified IMPES-type sequential scheme \cite{AzizSettari1979,CHEN2021113035,CHEN2019641} in which pressure and velocity are treated implicitly, concentration is advanced explicitly, and temperatures are solved implicitly. 
To stabilize
the advective terms in the transport equation for the concentration and enforce a discrete maximum principle, we incorporate a flux-corrected transport (FCT) strategy \cite{Kuzmin2012,Kuzmin-EG,kuzmin2023}. 
Time integration of the concentration equation is performed using a strong-stability-preserving second-order Runge–Kutta (SSP-RK2) method to achieve second-order temporal accuracy while preserving appropriate global and local bounds.

{We demonstrate that the proposed three-way LTNE model provides enhanced resolution of injected-fluid thermal evolution in fractured porous media relevant to EGS applications.} 
The numerical experiments illustrate how explicitly modeling injected-fluid heating reveals thermal structures that remain hidden in classical LTE or averaged LTNE formulations, thereby offering a more informative predictive tool for geothermal reservoir design and optimization.

The remainder of the paper is organized as follows. 
In Section~\ref{sec:model}, we present the three-way LTNE model and its coupling structure. 
Section~\ref{sec:numerics} describes the enriched Galerkin discretization, sequential splitting strategy, FCT stabilization, and SSP-RK2 \cite{GottliebShuTadmor2001} time integration. 
Numerical experiments for enhanced geothermal scenarios in fractured porous media are reported in Section~\ref{sec:numericsim}.

\section{Mathematical Modeling of Flow, Transport, and Heat Transfer}
\label{sec:model}
In this section, we describe the governing equations describing fluid pressure, and the local temperatures corresponding to the injected fluid, resident fluid, and solid matrix, where the injected fluid and resident fluid temperature are split by a concentration.
The coupled processes include Darcy flow, heat transfer between these interacting phases, and the species transport,



\subsection{Pressure equation}
The single-phase fluid flow in a porous medium $\mathbf{x}\in\Omega\subset\mathbb{R}^d$, with $d=2,3$, is governed by Darcy's law
\begin{equation}
\mathbf{u} = -\bm{\kappa}\nabla p,
\label{eq:Darcy-velocity}
\end{equation}
together with the mass conservation equation
\begin{equation}
\nabla\cdot\mathbf{u} := q_p \quad \text{in } \Omega,
\label{eq:mass-balance}
\end{equation}
where $\mathbf{u}:=\mathbf{u}(\mathbf{x})$ {is a $d$-dimensional vector field} denoting the Darcy velocity and 
$p:=p(\mathbf{x})$  represents the scalar fluid pressure. 
The term $q_p=q_p(\mathbf{x})$ denotes a source or sink. 
The mobility tensor is defined as $\bm{\kappa}=\mathbf{K}/\mu_f$, where 
$\mathbf{K}=K\mathbf{I}$ is the intrinsic permeability tensor, 
$K := K(\mathbf{x}) > 0$ is the scalar permeability, 
$\mu_f > 0$ is the fluid viscosity, and $\mathbf{I}$ is the $d\times d$ identity tensor.


The boundary conditions are prescribed as
\begin{equation}
p = p_D \quad \text{on } \partial \Omega_D,
\qquad
\mathbf{u} \cdot \mathbf{n} = g_N \quad \text{on } \partial \Omega_N,
\end{equation}
where $\partial \Omega$ {denotes} the boundary of $\Omega$, 
$\partial \Omega_D \cup \partial \Omega_N = \partial \Omega$ 
{and}
$\partial \Omega_D \cap \partial \Omega_N = \emptyset$.
Here, $\partial \Omega_D$ and $\partial \Omega_N$ denote the Dirichlet and Neumann 
{parts of the boundary}, with prescribed data 
$p_D := p_D(\mathbf{x})$ and $g_N := g_N(\mathbf{x})$, respectively.
{The vector $\mathbf{n}:=\mathbf n(\mathbf x)$ denotes the outward unit normal vector on $\partial \Omega$.}

\subsection{Temperature equation}\label{sec:Temperature}
Heat transfer in porous geothermal reservoir is often modeled using either the local thermal equilibrium (LTE) ~\cite{nield1992convection} assumption or the more general local thermal non-equilibrium (LTNE) formulation. The choice between these two models depends on the physical properties of the medium, flow regime, and the characteristic time scale of heat exchange between phases.

{
\subsubsection{Local thermal equilibrium (LTE) }\label{sec:LTE}
Under the local thermal equilibrium (LTE) assumption, the fluid temperature $T_f := T_f(\mathbf{x},t)$ and the solid temperature $T_s:= T_s(\mathbf{x},t)$ are equal and share a single temperature field $T:=T(\mathbf{x},t)$, i.e.,
$$
T_f = T_s = T.
$$
This assumption is valid when interphase heat exchange is sufficiently fast such that no appreciable temperature difference exists at the scale of a representative elementary volume (REV) ~\cite{kaviany1995,nield1992convection,Diersch2014}.

The governing LTE equation in a porous medium is written as 
\begin{equation}\label{eq:LTE}
    (\rho c)_{\mathrm{eff}} \frac{\partial T}{\partial t}
    + \rho_f c_f\, \mathbf{u} \cdot \nabla T
    - \nabla \cdot \left( \lambda_{\mathrm{eff}} \nabla T \right)
    = q_T
    \qquad \text{in } \Omega \times (0,t_{\text{final}}],
\end{equation}
with the initial condition
\begin{equation}
    T(\mathbf{x},0) = T^0(\mathbf{x})
    \qquad \text{in } \Omega,
\end{equation}
and the Dirichlet and Neumann boundary conditions
\begin{equation}
    T = T_D 
    \quad \text{on } \Gamma_D \times (0,t_{\text{final}}],
    \qquad
    -\lambda_{\mathrm{eff}} \nabla T \cdot \mathbf{n} = g_T
    \quad \text{on } \Gamma_N \times (0,t_{\text{final}}].
\end{equation}

Here, $t_\text{final}$ is final computation time, $q_T := q_T(\mathbf{x},t) \in \mathbb{R}$ denotes a source or sink term. The quantity 
$(\rho c)_{\mathrm{eff}}$ represents the effective volumetric heat capacity of the porous medium, 
obtained by volume-averaging the contributions of the fluid and solid phases,
\begin{equation}
(\rho c)_{\mathrm{eff}}
:= \phi\, \rho_f c_f
+ (1-\phi)\, \rho_s c_s,
\end{equation}
where $0 \le \phi \le 1$ is the porosity, $\rho_f>0$ and $\rho_s>0$ are the densities, and 
$c_f>0$ and $c_s>0$ are the specific heat capacities of the fluid and solid, respectively.

The term $\lambda_{\mathrm{eff}}$ denotes the effective thermal conductivity of the porous medium,
defined as the sum of the fluid and solid contributions,
\begin{equation}
\lambda_{\mathrm{eff}}
:= \lambda_{f,\mathrm{eff}}
+ \lambda_{s,\mathrm{eff}},
\end{equation}
where $\lambda_{f,\mathrm{eff}} := \phi \lambda_f$ and 
$\lambda_{s,\mathrm{eff}} := (1-\phi)\lambda_s$ are the effective thermal conductivities of the fluid and solid phases, based on their intrinsic thermal conductivities $\lambda_f >0 $ and $\lambda_s>0$.

\subsection{Local thermal non-equilibrium (LTNE)}

In contrast to the LTE assumption, the LTNE framework allows the fluid and solid phases to have distinct temperature fields. 
This distinction becomes {important} in geothermal and subsurface heat-transport problems where cold fluid is injected into a much hotter formation. 
In such situations, interphase heat transfer is not instantaneous, and significant temperature differences may exist at the {scale of a representative elementary volume (REV)}. 
To account for this finite-rate {interphase heat transfer}, the LTNE model introduces two separate temperature fields: the fluid temperature $T_f$ and the solid temperature $T_s$.


The governing equations for the fluid and solid temperatures $T_f$ and $T_s$ are {given by}
\begin{align}
    \phi \rho_f c_f \dfrac{\partial T_f}{\partial t} + \rho_f c_f \mathbf u_f \cdot \nabla T_f - \nabla \cdot (\lambda_{f,\eff} \nabla T_f) &= h_{fs}a_{fs}(T_s - T_f) + q_f\label{eq:LTNE-f}  &&\text{in } \Omega \times (0,t_{\text{final}}],\\
    (1-\phi) \rho_s c_s \dfrac{\partial T_s}{\partial t} - \nabla \cdot \left(\lambda_{s,\eff} \nabla T_s\right) &= h_{fs}a_{fs}(T_f - T_s) + q_s\label{eq:LTNE-s}  &&\text{in } \Omega \times (0,t_{\text{final}}],
\end{align}
where $q_f:=q_f(\mathbf x,t)$, $q_s:=q_s(\mathbf x,t)$ are the volumetric heat source or sink term associated with fluid and solid phases, and $h_{fs}>0$ is the interfacial heat transfer coefficient and $a_{fs}>0$ denotes the specific interfacial area {(i.e., the solid--fluid surface area per unit bulk volume)}. 
The term $h_{fs}a_{fs}(T_s - T_f)$ represents {the finite-rate heat exchange} between the fluid and solid phases. 
In particular, when the {interphase heat transfer coefficient} $h_{fs} a_{fs}$ is sufficiently large, the two temperature fields rapidly approach thermal equilibrium, and the LTNE model {reduces to} the LTE formulation.


The initial temperatures for the fluid and solid are prescribed to be identical:
\begin{equation}
    T_f(\mathbf x,0) = T_s(\mathbf x, 0) = T^0(\mathbf x) \quad \text{in } \Omega.
\end{equation}
For the fluid phase, a Dirichlet boundary condition is imposed on the inflow boundary, while a Neumann condition is prescribed elsewhere:
\begin{equation}
    T_f = T_\text{in} 
    \quad \text{on } \Gamma_\text{in} \times (0,t_\text{final}],\qquad 
    -\lambda_{f.\eff} \nabla T_f \cdot \mathbf n = g_f 
    \quad \text{on } \partial \Omega \times (0,t_\text{final}].
\end{equation}
For the solid phase, only Neumann boundary conditions are considered:
\begin{equation}
    -\lambda_{s.\eff} \nabla T_s \cdot \mathbf n = g_s 
    \quad \text{on } \partial \Omega \times (0,t_\text{final}].
\end{equation}

{In the next section, we present an asymptotic analysis demonstrating how the LTNE model reduces to the LTE model in the limit of strong interfacial heat exchange.}

\subsubsection{Asymptotic limit of strong fluid--solid heat exchange}\label{subsec:asymptotic}
Define the phase heat capacities per unit REV
\[
C_f := \phi\,\rho_f c_f, \qquad C_s := (1-\phi)\,\rho_s c_s,
\]
the interfacial coupling \(\Lambda := h_{fs} a_{fs}\), and the
{mixture temperature} $T_\text{mix} := T_\text{mix}(\mathbf x,t)$ as the energy-equivalent average
\begin{equation}
(C_f + C_s)\,T_{\mathrm{mix}} \;:=\; C_f T_f + C_s T_s.
\label{eq:Tmix-def}
\end{equation}
Summing \eqref{eq:LTNE-f}--\eqref{eq:LTNE-s} and using \eqref{eq:Tmix-def} yields the exact mixture equation
\begin{equation}
(C_f + C_s)\,\partial_t T_{\mathrm{mix}} \;+\; \rho_f c_f \,\mathbf{u}_f\!\cdot\!\nabla T_f
\;=\; \nabla\!\cdot\!\big(\lambda_{f,\eff}\nabla T_f
+ \lambda_{s,\eff} \nabla T_s\big) \;+\; q_f + q_s,
\label{eq:Tmix-eq}
\end{equation}
where $\lambda_{f,\eff} = \phi\lambda_f$ and $\lambda_{s,\eff}=(1-\phi)\lambda_s$, as before. 

We introduce the temperature difference \(\theta := T_f - T_s\) and use \eqref{eq:Tmix-def} to write
\begin{equation}
T_f = T_{\mathrm{mix}} + \frac{C_s}{C_f+C_s}\,\theta, 
\qquad
T_s = T_{\mathrm{mix}} - \frac{C_f}{C_f+C_s}\,\theta.
\label{eq:T-split}
\end{equation}
Then subtracting \eqref{eq:LTNE-s} divided by $C_s$ from \eqref{eq:LTNE-f} divided by $C_f$, we obtain the evolution equation 
\begin{equation}
\partial_t \theta
\;+\; \dfrac{1}{\phi}\mathbf{u}_f\!\cdot\!\nabla T_f
\;-\; \frac{1}{C_f}\,\nabla\!\cdot\!\big(\lambda_{f,\eff}\nabla T_f\big)
\;+\; \frac{1}{C_s}\,\nabla\!\cdot\!\big(\lambda_{s,\eff}\nabla T_s\big)
\;+\; \Big(\frac{q_f}{C_f}-\frac{q_s}{C_s}\Big)
\;=\; -\,\gamma\,\theta
\label{eq:theta}
\end{equation}
for $\theta$ with the decay rate
\begin{equation}
\gamma \;:=\; \Lambda\!\left(\frac{1}{C_f}+\frac{1}{C_s}\right)
= h_{fs} a_{fs}\!\left(\frac{1}{\phi\,\rho_f c_f}
+\frac{1}{(1-\phi)\,\rho_s c_s}\right).
\label{eq:gamma}
\end{equation}

Let \(\varepsilon := \gamma^{-1}\) be small, corresponding to strong interfacial exchange \(\Lambda\to\infty\).
Assuming coefficients and sources remain \(O(1)\) as \(\varepsilon\to 0\), equation \eqref{eq:theta} implies
\(\theta = O(\varepsilon)\) after a fast transient \(t=O(\varepsilon)\).
Then it follows from \eqref{eq:T-split} that
\[
T_f = T_{\mathrm{mix}} + O(\varepsilon), \qquad
T_s = T_{\mathrm{mix}} + O(\varepsilon).
\]
Invoking \eqref{eq:Tmix-eq} and taking \(\varepsilon\to 0\) gives 
\begin{equation}
(C_f + C_s)\,\partial_t T_{\text{eq}} \;+\; \rho_f c_f \,\mathbf{u}\!\cdot\!\nabla T_\text{eq}
\;=\; \nabla\!\cdot\!\Big(\,[\lambda_f^{\mathrm{eff}}+\lambda_s^{\mathrm{eff}}]\,\nabla T_{\text{eq}}\Big)
\;+\; q_f + q_s,
\label{eq:LTE-limit}
\end{equation}
where \(T_{\text{eq}} := \lim_{\varepsilon\to 0} T_{\mathrm{mix}}\).
Since $C_f+C_s = (\rho c)_\eff$, the limit \eqref{eq:LTE-limit} matches the LTE equation \eqref{eq:LTE}.

As \(h_{fs}a_{fs}\!\to\!\infty\), the temperature difference \(\theta=T_f-T_s\) relaxes as
\(\theta=O\!\left((h_{fs}a_{fs})^{-1}\right)\), and the mixture temperature
\(T_{\mathrm{mix}}\) converges to the single LTE temperature \(T_{\text{eq}}\) solving \eqref{eq:LTE-limit}.

}

{
\subsection{Three-way LTNE coupling}\label{sec:three-way-LTNE}

To accurately capture the thermal behavior of geothermal systems with strong thermal contrasts—such as cold-fluid injection into a hot reservoir—it is often necessary to distinguish not only between the fluid and solid phases, but also between {the injected and resident fluid components}. 
These two fluid populations may coexist within the same pore space while exhibiting distinct temperatures due to differences in arrival times, velocities, and local thermal exposure. 
Traditional LTE or {two-temperature LTNE} models are not able to explicitly represent such thermal separation. 
Therefore, we introduce a {three-temperature LTNE framework} with temperature fields for the injected fluid $T_i(\mathbf{x},t)$, the resident fluid $T_r(\mathbf{x},t)$, and the solid matrix $T_s(\mathbf{x},t)$, together with a concentration variable $0 \leq z \leq 1$ representing the injected-fluid fraction.

Here, {the volumetric enthalpy of the fluid phase is decomposed as}
\begin{equation}\label{eq:enthalpy for fluid}
    C_f T_f = C_i T_i + C_r T_r,
\end{equation}
where the volumetric heat capacities of the injected and resident fluid components are defined by
\begin{equation}
    C_i := (\phi z)\,\rho_f c_f,
    \qquad 
    C_r := \big(\phi(1-z)\big)\,\rho_f c_f.
\end{equation}
The effective thermal conductivities of the injected and resident fluid components are defined as
\begin{equation}
    \lambda_{i,\mathrm{eff}} := \phi z\,\lambda_f,
    \qquad 
    \lambda_{r,\mathrm{eff}} := \phi(1-z)\,\lambda_f,
\end{equation}
so that
\begin{equation*}
    \lambda_{i,\mathrm{eff}} + \lambda_{r,\mathrm{eff}}
    = \phi\lambda_f = \lambda_{f,\mathrm{eff}}.
\end{equation*}

The fluid--solid interfacial heat exchange is partitioned according to the injected-fluid concentration $z$:
\begin{equation}
    h_{fs}a_{fs} = (z\,h_{fs}a_{fs}) + \big((1-z)\,h_{fs}a_{fs}\big)
    = \Lambda_{is} + \Lambda_{rs},
\end{equation}
where 
\begin{equation}
    \Lambda_{is} := z\,h_{fs}a_{fs},
    \qquad
    \Lambda_{rs} := (1-z)\,h_{fs}a_{fs}
\end{equation}
represent the injected--fluid--solid and resident--fluid--solid coupling strengths, respectively. 

In addition, we account for direct heat exchange between the injected and resident fluid components.
The strength of this interaction is assumed to scale with both the relative proportions of the two fluids and the available pore volume 
{within which} they coexist.
Therefore, we define
\begin{equation}
    \Lambda_{ir} := \phi z(1-z) h_{ir} a_{ir},
\end{equation}
so that the injected--resident heat exchange scales with the product 
$z(1-z)$, {which vanishes when either component is absent,}
and is further modulated by the porosity $\phi$, 
{representing}
the total pore space available for fluid--fluid thermal interaction.

With these definitions, the three-temperature LTNE system for the injected fluid, resident fluid, and solid matrix is written as
\begin{align}
    C_i \frac{\partial T_i}{\partial t} 
    +z\rho_f c_f  \mathbf{u}\cdot\nabla T_i 
    - \nabla\cdot\!\left(\lambda_{i,\mathrm{eff}}\,\nabla T_i\right)
    &= \Lambda_{is}(T_s - T_i)+\Lambda_{ir}(T_r - T_i) + q_i
    &&{\text{in } \Omega \times (0,t_{\text{final}}]},\label{eq:Ti}\\
    C_r \frac{\partial T_r}{\partial t} 
    +(1-z)\rho_f c_f  \mathbf{u}\cdot\nabla T_r 
    - \nabla\cdot\!\left(\lambda_{r,\mathrm{eff}}\,\nabla T_r\right)
    &= \Lambda_{rs}(T_s - T_r)+\Lambda_{ir}(T_i - T_r) + q_r
    &&{\text{in } \Omega \times (0,t_{\text{final}}]},\label{eq:Tr}\\
    C_s\,\frac{\partial T_s}{\partial t}
    - \nabla\cdot\!\left(\lambda_{s,\mathrm{eff}} \nabla T_s\right)
    &= \Lambda_{is}(T_i - T_s) + \Lambda_{rs}(T_r - T_s) + q_s
    &&{\text{in } \Omega \times (0,t_{\text{final}}]}. \label{eq:Ts}
\end{align}
Here, $q_i := q_i(\mathbf{x},t)$, $q_r := q_r(\mathbf{x},t)$, and $q_s := q_s(\mathbf{x},t)$
denote volumetric heat source or sink terms associated with the injected fluid,
resident fluid, and solid phases, respectively.

For compactness and thermodynamic consistency, the initial states are prescribed 
in terms of the corresponding volumetric enthalpies:
\begin{equation}
T_i(\mathbf x,0) = \frac{H_i^0(\mathbf x)}{C_i},
\qquad
T_r(\mathbf x,0) = \frac{H_r^0(\mathbf x)}{C_r},
\qquad
T_s(\mathbf x,0) = \frac{H_s^0(\mathbf x)}{C_s},
\qquad
{\text{in } \Omega,}
\end{equation}
where 
\[
H_i = C_i T_i, 
\qquad 
H_r = C_r T_r, 
\qquad 
H_s = C_s T_s
\]
denote the volumetric enthalpies of the injected fluid, resident fluid, and solid phases, respectively. 
In particular, we set $H_i^0(\mathbf x) = 0$, so that the injected-fluid subsystem {initially contains no thermal energy}.

Let $
\Gamma_{\mathrm{in}} := \{ \mathbf{x} \in \partial \Omega : \mathbf{u}_f(\mathbf{x}) \cdot \mathbf{n}(\mathbf{x}) < 0 \}$.
Injection at a given temperature \(T_{\mathrm{inj}}:=T_\mathrm{inj}(\mathbf x)\) is imposed equivalently by the enthalpy boundary condition
\begin{equation}
T_{\mathrm{inj}} = H_{i,\mathrm{inj}}/C_i 
\quad\text{on } \Gamma_{\mathrm{in}}
\end{equation}
On $\Gamma_{\mathrm{out}} := \partial \Omega \setminus \Gamma_{\mathrm{in}}$, we
impose the Neumann boundary condition
\begin{equation}
    -\lambda_{i,\eff}\nabla T_i = g_i\quad  \text{on $\Gamma_\text{out}$}.
\end{equation} 
For solid and resident fluid temperatures, we impose full Neumann boundary conditions
\begin{equation}
   - \lambda_{r,\eff} \nabla T_r = g_r,\quad -\lambda_{s,\eff}\nabla T_s = g_s\quad \text{on }\partial \Omega.
\end{equation}

In the next section, we formulate an evolution equation for the concentration of injected fluid, which is used in the three-way-LTNE modeling.

\subsection{Concentration equation}
 The injected fluid contains a dissolved or miscible component whose concentration is described by the volume fraction
\[
z := z(\mathbf{x}, t).
\]
We consider a binary miscible mixture, so the volume fractions of the two components satisfy
\[
z + (1 - z) = 1,
\]
meaning that the total fluid volume fraction is unity.

The evolution of \(z\) is governed by the advection equation 
\begin{equation}
\frac{\partial (\phi \rho_f z)}{\partial t}
+ \nabla \cdot (\rho_f z\, \mathbf{u})
= 0\quad{\text{in } \Omega \times (0, t_{\mathrm{final}}]},\label{eq:z}
\end{equation}
where $\phi$ is the porosity, $\rho_f$ is the fluid density,
and $q_z:=q_z(\mathbf x,t)$ represents volumetric source or sink terms.

The initial condition for the species concentration is
\begin{equation}
    z(\mathbf x,0) = z^0(\mathbf x) \quad \text{in } \Omega.
\end{equation}
We impose the inflow boundary condition 
\begin{equation}
z = z_{\mathrm{in}}
\quad \text{on } \Gamma_{\mathrm{in}} \times (0, t_{\mathrm{final}}],
\end{equation}
where
$
\Gamma_{\mathrm{in}} := \{ \mathbf{x} \in \partial \Omega : \mathbf{u}_f(\mathbf{x}) \cdot \mathbf{n}(\mathbf{x}) < 0 \}$, as before.

\section{Numerical Methods} \label{sec:numerics}

In this section, we present the numerical discretization of the governing equations. 
Spatial discretization is performed using the enriched Galerkin (EG) finite element method, which combines continuous $\mathbb{Q}_1$ elements with a discontinuous $\mathbb{Q}_0$ enrichment to ensure local conservation \cite{becker2003,LEE2016}. 
Implicit and explicit time integration schemes are employed depending on the stability requirements. 
An EG formulation is used for the pressure, temperature, and concentration equations, and a flux-corrected transport (FCT) scheme \cite{Kuzmin-EG} is applied to stabilize
the advective terms in the transport equation for the concentration.
\subsection{Enriched Galerkin method}
Let $\Tau_h = \{K_e\}_{e=1}^{E_h}$ denote a non-degenerate partition of the domain $\Omega$ into $E_h$ rectangular cells of maximum diameter $h=\max_{1\leq e\leq E_h} h_{K_e}$, where $h_{K_e}$ is the diameter of $K_e$. The vertices of $\Tau_h$ are denoted by $\mathbf x_1,\ldots,\mathbf x_{N_h}$. We store the indices of vertices belonging to a given cell $K_e$ in the integer set $\mathcal N^e$ and the indices of elements that contain a given vertex $\mathbf x_i$ in the integer set $\mathcal E_i$. The set $\mathcal N_i=\bigcup_{e\in\mathcal E_i}\mathcal N^e$ contains the indices of all vertices belonging to at least one cell that contains~$\mathbf x_i$. Note that $i\in\mathcal N_i$. The boundary of $K_e$ is denoted by $\partial K_e$ with outward normal vector denoted by~$\mathbf n_e$.

The set of all edges $\gamma$ in the collection of elements is denoted by $\Gamma_h$, and the length of $\gamma$ is denoted by $h_\gamma$. The sets of interior, Dirichlet boundary, Neumann boundary, and inflow/outflow boundary edges are denoted by $\Gamma_h^\text{I},\Gamma_h^D,\Gamma_h^N,\Gamma_h^\text{in}$, and $\Gamma_h^\text{out}$, respectively. Each edge is associated with a normal vector $\mathbf n_\gamma$. For a boundary edge, $\mathbf n_\gamma$
represents the outward normal to $\partial\Omega$.

Let $\mathbb Q_{1}(\hat K)$ denote the space of multilinear polynomials $\hat v:\hat K\to\R$ defined on the reference element $\hat K=[0,1]^d$. Using a multilinear mapping $F_e:\hat K\to K_e\in\mathcal T_h$, we construct the space $\mathbb Q_{1}(K_e)$ of polynomials $v:K_e\to\R$ such that $v=\hat v\circ F_e^{-1}$ for some $\hat v\in\mathbb Q_{1}(\hat K)$. 

The finite element space CG-$\mathbb Q_1$ of the classical continuous Galerkin method using $\mathbb Q_1$ Lagrange elements on the partition $\Tau_h$ is defined as
\begin{equation}
V^{\text{CG}}_h := \{ v_h \in L^2(\Omega)\,:\, v_h|_{K_e} \in \mathbb Q_1(K_e)\ \forall K_e \in \Tau_h \} \cap \mathbb C(\bar \Omega),
\end{equation}
where $\mathbb C(\bar\Omega)$ denotes the space of functions that are
continuous on $\bar \Omega$. The DG-$\mathbb Q_0$ space

\begin{equation}
  V^{\DG}_{h} := \{\delta v_h\in L^2(\Omega)\,:\,
  \delta v_h|_{K_e}\in\mathbb Q_0(K_e)\ \forall K_e\in \Tau_h\}
\end{equation}
consists of functions that are constant on elements of the partition $\Tau_h$. 
The finite element space of the EG-$\mathbb Q_1$ method is then defined as  (cf. \cite{becker2003,Kuzmin-EG,LEE2016,LEE201719})
\begin{equation}
    V_{h} := V^{\text{CG}}_h\oplus V^{\DG}_{h}.
\end{equation}

Finally, we introduce the jump and average values for $v_h\in V_h$ on edges $\gamma \in \Gamma_h$. For an interior edge $\gamma\in\Gamma_h^\text{I}$ such that $\gamma=\partial K_e\cap \partial K_{e'}$ with the associated normal vector $\mathbf n_\gamma = \mathbf n_e$, the jump value is defined as
\begin{equation}
    \jump{v_h} := 
        \left(v_h|_{K_e}\right)\Big|_\gamma - \left(v_h|_{K_{e'}}\right)\Big|_\gamma.
\end{equation}
Similarly, the average value $\av{\varphi}$ is defined as 
\begin{equation}
    \av{v_h} := 
        \frac{1}{2}\left(v_h|_{K_e}\right)\Big|_\gamma +\frac{1}{2} \left(v_h|_{K_{e'}}\right)\Big|_\gamma.
\end{equation}
For edges on the boundary $\gamma \subset \partial \Omega$, we define $\jump{v_h} = \av{v_h} = v_h|_\gamma$.

\subsection{Temporal discretization}\label{sec:time-discretization}
The time discretization is carried out on the time interval
$[0,t_{\mathrm{final}}]$, with $t_{\mathrm{final}}>0$, using a given number of
uniform time steps $N\in\mathbb{N}$. The time step size is defined as
$\Delta t = t_{\mathrm{final}}/N$, and the discrete time levels are given by
$t^n := n\Delta t$ for $0 \le n \le N$.
A time-dependent function evaluated at time $t^n$ is denoted by
$\phi^n := \phi(\cdot,t^n)$.

Let $\phi(t)$ satisfy the evolution equation
$\partial_t \phi = \mathcal{F}(\phi)$.
If an implicit time discretization is preferable for stability reasons,
we use the first-order backward Euler method
\begin{equation}\label{scheme:be}
\frac{\phi^{n+1}-\phi^n}{\Delta t}
=
\mathcal{F}(\phi^{n+1}),
\qquad n \ge 0
\end{equation}
or the second-order backward differentiation formula (BDF2) 
\begin{equation}\label{scheme:bdf2}
\frac{3\phi^{n+1}-4\phi^n+\phi^{n-1}}{2\Delta t}
=
\mathcal{F}(\phi^{n+1}),
\qquad n \ge 1.
\end{equation}

In addition to implicit schemes, we consider an explicit strong-stability-preserving (SSP) Runge--Kutta method of second order (RK2). It performs two forward Euler steps
\begin{equation}\label{scheme:rk2-1stage}
\begin{cases}
    \phi^{(1)}
=
\phi^n
+
\Delta t\,\mathcal{F}(\phi^n),\\
\phi^{(2)}
=
\phi^{(1)}
+
\Delta t\,\mathcal{F}(\phi^{(1)})
\end{cases}
\end{equation}
and updates the solution as follows:
\begin{equation}\label{scheme:rk2-update}
\phi^{n+1}
=
\dfrac{1}{2}(\phi^{(1)}+\phi^{(2)}).
\end{equation}

\subsection{Spatial discretization for pressure}
The EG finite element approximation to the pressure $p$ that defines the Darcy velocity \eqref{eq:Darcy-velocity} is denoted by $p_h := p_h(\mathbf x)\in V_h$. The spatial semi-discretization of the Poisson equation for $p$, which follows from \eqref{eq:Darcy-velocity} and \eqref{eq:mass-balance}, employs the EG weak form
\begin{equation}\label{FEM-p}
    \text{find } p_h \in V_{h} \text{ such that }\mathcal A_p(p_h, v_h) = \mathcal F_p(v_h),\qquad \forall v_h \in V_{h},
\end{equation}
where $\mathcal A_p$ and $\mathcal F_p$ are the bilinear form and linear functional defined as 
\begin{equation}\label{eq:pressure-bilinear}
\begin{split}
        \mathcal A_{p}(p_h,v_h) =& 
        \sum_{e=1}^{E_h} \int_{K_e} \bm \kappa\nabla p_h\cdot\nabla v_h\,\mathrm d\mathbf x 
        + \sum_{\gamma\in \Gamma_h^\mathrm I \cup \Gamma_h^\mathrm D}\dfrac{\alpha_p}{h_\gamma}\int_{\gamma}\jump{\bm  \kappa p_h}\jump{v_h}\mathrm ds\\
       & - \sum_{\gamma\in \Gamma_h^\mathrm I\cup\Gamma_h^\mathrm D } \int_\gamma \av{\bm \kappa \nabla p_h\cdot \mathbf n_\gamma}\jump{v_h} \mathrm ds + \theta \sum_{e\in\Gamma_h^\mathrm I \cup\Gamma_h^\mathrm D}  \int_\gamma \av{\bm \kappa \nabla v_h\cdot \mathbf n_\gamma} \jump{p_h} \mathrm ds,
\end{split}
\end{equation}
and 
\begin{equation}\label{eq:pressure-linearFunctional}
\begin{split}
    \mathcal F_p(v_h) =&
    \sum_{e=1}^{E_h}\int_{K_e} q_p v_h \,\mathrm d\mathbf x +  \sum_{\gamma\in \Gamma_h^\mathrm D}\dfrac{\alpha_p}{h_\gamma}\int p_{\!D} \jump{v_h}\,\mathrm ds 
    - \sum_{\gamma\in\Gamma_h^\mathrm N}\int_{\gamma}g_N\,\jump{v_h} \mathrm ds + \theta \sum_{\gamma\in \Gamma_h^\mathrm D}\int_\gamma p_{\!D} \av{\bm \kappa\nabla v_h\cdot \mathbf n_\gamma}\,\mathrm ds .
\end{split}
\end{equation}
Here $\alpha_p>0$ is a penalty parameter. The choice of $\theta$ leads to symmetric interior penalty Galerkin (SIPG) with $\theta = -1$, incomplete interior penalty Galerkin (IIPG) with $\theta = 0$, and non-symmetric interior penalty Galerkin (NIPG) with $\theta = 1$. Throughout this paper, we choose IIPG $(\theta = 0)$. 


\subsubsection{Locally conservative Darcy's velocity}
The EG approximation to the Darcy velocity $\mathbf u_f$ is denoted by $\mathbf u_h$. 
We define a globally and locally conservative $\mathbf u_h$ by differentiating the EG pressure $p_h$; cf.~\cite{sun2005discontinuous,LEE201719}
\begin{subequations}
\begin{alignat}{2}
    \mathbf u_h|_{K_e} &= - \bm \kappa \nabla p_h & \quad \forall K_e\in\Tau_h,\\
    (\mathbf u_h \cdot \mathbf n_\gamma)|_\gamma&= -\av{\bm \kappa \nabla p_h\cdot \mathbf n_\gamma} + \dfrac{\alpha_p}{h_\gamma} \jump{\bm \kappa p_h} \quad&\forall \gamma \in \Gamma_h^\mathrm I,\\
    (\mathbf u_h \cdot \mathbf n_\gamma)|_\gamma&= g_N \quad &\forall \gamma\in \Gamma_h^\mathrm N,\\
    (\mathbf u_h \cdot \mathbf n_\gamma)|_\gamma&= -{\bm \kappa \nabla p_h\cdot \mathbf n_\gamma} + \dfrac{\alpha_p}{h_\gamma} \bm \kappa (p_h-p_{\!D}) \quad&\forall \gamma \in \Gamma_h^\mathrm D.
\end{alignat}
\end{subequations}
This definition guarantees local mass conservation on each element $K_e \in \mathcal T_h$ in the sense that 
\begin{equation}
    \sum_{\gamma\in\partial K_e}\int_{\gamma} \mathbf u_h \cdot \mathbf n_\gamma\, \mathrm ds = \int_{K_e} q_m\mathrm d\mathbf x.
\end{equation}

\subsection{Spatial discretization for three-way-LTNE}
We denote the EG solution to LTNE system \eqref{eq:Ti}--\eqref{eq:Ts} by $T_{i,h},T_{r,h},T_{s,h}$. The spatial semi-discretization of the full LTNE system is similar to \eqref{eq:FEM-Concentration} and reads:
$$
\text{find }(T_{i,h},T_{r,h},T_{s,h})\in V_h\times V_h\times V_h \ \text{such that } 
$$
\begin{align}
    \sum_{e=1}^{E_h} \int_{K_e}\phi_{i}\rho_{i}c_i z\dfrac{\partial T_{i,h}}{\partial t} v_h\,\mathrm dx &=
    \mathcal A^\text{adv}_i(T_{i,h},v_h;z) 
    + \mathcal A^\text{diffu}_i(T_{i,h},v_h;z)
    +H_{si}(v_h;z)+H_{ri}(v_h;z),\\
    \sum_{e=1}^{E_h} \int_{K_e}\phi_{r}\rho_{r}c_r (1-z)\dfrac{\partial T_{r,h}}{\partial t} w_h\,\mathrm dx &=
    \mathcal A^\text{adv}_r(T_{r,h},w_h;z) 
    + \mathcal A^\text{diffu}_r(T_{r,h},w_h;z)
    +H_{sr}(w_h;1-z)+H_{ir}(w_h;1-z), \\
    \sum_{e=1}^{E_h} \int_{K_e}(1-\phi)\rho_{s}c_s \dfrac{\partial T_{s,h}}{\partial t} s_h\,\mathrm dx &=
     \mathcal A^\text{diffu}_s(T_{s,h},s_h;z)
     +H_{is}(s_h;z)+H_{rs}(s_h;1-z)
\end{align}
 for all $(v_h,w_h,s_h)\in V_h\times V_h\times V_h$. The advective contributions to the semi-discrete equations for the injected and resident fluid are denoted by $\mathcal A_i^\text{adv}$ and $\mathcal A_r^\text{adv}$, respectively. These functionals depend on the concentration $z$ and are defined as
\begin{equation}
\begin{split}
    \mathcal A_i^\text{adv}( T_{i,h}, v_h;z):=&  \sum_{e=1}^{E_h}\int_{K_e} \phi_i(\rho c)_i z\,T_{i,h} \mathbf u_h \cdot \nabla v_h\,\mathrm d\mathbf x - \sum_{e=1}^{E_h} \int_{\partial K_e^+} \phi_i(\rho c)_i z\,T_{i,h}^+ \mathbf u_h \cdot \mathbf n_e\,v_h\,\mathrm ds  \\
    &- \sum_{\gamma \in \Gamma_h^\text{in}}\int_\gamma \phi_i(\rho c)_i z\,T_{\text{in}} v_h\,\mathrm ds,
\end{split}
\end{equation}
\begin{equation}
\begin{split}
    \mathcal A_r^\text{adv}( T_{r,h}, w_h;z):=&  \sum_{e=1}^{E_h}\int_{K_e} \phi_r(\rho c)_r (1-z)T_{r,h} \mathbf u_h \cdot \nabla w_h\,\mathrm d\mathbf x - \sum_{e=1}^{E_h} \int_{\partial K_e^+} \phi_r(\rho c)_r (1-z)T_{r,h}^+ \mathbf u_h \cdot \mathbf n_e\,w_h\,\mathrm ds  \\
    &- \sum_{\gamma \in \Gamma_h^\text{in}}\int_\gamma \phi_r(\rho c)_r (1-z)T_{\text{in}} w_h\,\mathrm ds.
\end{split}
\end{equation}
The diffusive contributions to the temperature equations for the three phases are represented by
\begin{equation}
\begin{split}
    \mathcal A_i^\text{diffu}( T_{i,h}, v_h;z):=& - \sum_{e=1}^{E_h} \int_{K_e} \lambda_{i,\eff} \nabla T_{i,h}\cdot \nabla v_h\,\mathrm dx\\
    &- \sum_{\gamma\in \Gamma_h^I} \dfrac{\alpha_i}{h_\gamma} \int_\gamma \jump{\lambda_\mathrm{i,eff} T_{i,h}}\jump{v_h}\,\mathrm ds
    + \sum_{\gamma \in \Gamma_h^\mathrm I} \av{\lambda_\mathrm{i,eff}\nabla T_{i,h}\cdot \mathbf n_\gamma} \jump{v_h}\,\mathrm ds \\
    &- \sum_{\gamma \in \partial \Omega} \int_\gamma \left(\lambda_\mathrm{i,eff}\,g_{T_{i}} - \av{\lambda_\mathrm{i,eff} \nabla T_{i,h}\cdot \mathbf n_\gamma}\right)\jump{v_h}\,\mathrm ds,
\end{split}
\end{equation}
\begin{equation}
\begin{split}
    \mathcal A_r^\text{diffu}( T_{r,h}, w_h;z):=& - \sum_{e=1}^{E_h} \int_{K_e} \lambda_{i,\eff} \nabla T_{r,h}\cdot \nabla w_h\,\mathrm dx\\
    &- \sum_{\gamma\in \Gamma_h^I} \dfrac{\alpha_r}{h_\gamma} \int_\gamma \jump{\lambda_\mathrm{r,eff} T_{r,h}}\jump{w_h}\,\mathrm ds
    + \sum_{\gamma \in \Gamma_h^\mathrm I} \av{\lambda_\mathrm{r,eff}\nabla T_{r,h}\cdot \mathbf n_\gamma} \jump{w_h}\,\mathrm ds \\
    &- \sum_{\gamma \in \partial \Omega} \int_\gamma \left(\lambda_\mathrm{r,eff}\,g_{T_{r}} - \av{\lambda_\mathrm{r,eff} \nabla T_{r,h}\cdot \mathbf n_\gamma}\right)\jump{w_h}\,\mathrm ds,
\end{split}
\end{equation}
and
\begin{equation}
\begin{split}
    \mathcal A_s^\text{diffu}( T_{s,h}, s_h;z):=& - \sum_{e=1}^{E_h} \int_{K_e} \lambda_{s,\eff} \nabla T_{i,h}\cdot \nabla s_h\,\mathrm dx\\
    &- \sum_{\gamma\in \Gamma_h^I} \dfrac{\alpha_s}{h_\gamma} \int_\gamma \jump{\lambda_\mathrm{s,eff} T_{s,h}}\jump{s_h}\,\mathrm ds
    + \sum_{\gamma \in \Gamma_h^\mathrm I} \av{\lambda_\mathrm{s,eff}\nabla T_{s,h}\cdot \mathbf n_\gamma} \jump{s_h}\,\mathrm ds \\
    &- \sum_{\gamma \in \partial \Omega} \int_\gamma \left(\lambda_\mathrm{s,eff}\,g_{T_{s}} - \av{\lambda_\mathrm{s,eff} \nabla T_{s,h}\cdot \mathbf n_\gamma}\right)\jump{s_h}\,\mathrm ds,
\end{split}
\end{equation}
where $\alpha_i,\alpha_r,\alpha_s >0$ are penalty parameters. 
The heat exchange between the phases is taken into account through 
\begin{align}
    H_{si}(v_h;z) &:= \sum_{e=1}^{E_h}\int_{K_e}zh_{is}a_{is}(T_s - T_i) v_h\,\mathrm dx,\\
    H_{ri}(v_h;1 - z) &:= \sum_{e=1}^{E_h}\int_{K_e}(1-z)h_{ir}a_{ir}(T_r - T_i) v_h\,\mathrm dx,\\
    H_{sr}(v_h;1-z) &:= \sum_{e=1}^{E_h}\int_{K_e}(1-z)h_{rs}a_{rs}(T_s - T_r) w_h\,\mathrm dx
\end{align}
with 
$$
H_{is}(\,\cdot\,;\,\cdot\,) = -H_{si}(\,\cdot\,;\,\cdot\,),\quad H_{ri} = (\,\cdot\,;\,\cdot\,) = -H_{ir}(\,\cdot\,;\,\cdot\,),\quad H_{rs} (\,\cdot\,;\,\cdot\,) = - H_{sr} (\,\cdot\,;\,\cdot\,).
$$

We use the implicit BDF2 scheme for time discretization in the temperature equations, because an explicit treatment of diffusive terms would result in severe time step restrictions.  

\subsection{Spatial discretization for concentration}
The EG finite element approximation to the solution of \eqref{eq:z} is denoted by $z_h=z_h(\mathbf x,t) \in V_h$. The spatial semi-discretization of \eqref{eq:z} is based on the following weak form:   
\begin{equation}\label{eq:FEM-Concentration}
\begin{split}
\mbox{find $z_h \in V_h$ such that} \qquad
    \sum_{e = 1}^{E_h}& \int_{K_e} \phi\rho_f\dfrac{\partial z_h}{\partial t} v_h\,\mathrm dx = \mathcal A_z^\text{adv}(z_h,v_h) 
    \qquad 
    \forall v_h\in V_h,
\end{split}
\end{equation}
where $\mathcal A_z^\text{adv}(z_h,v_h):=\mathcal A_z^\text{adv}(z_h,v_h;\mathbf u_h)$ 
and it is defined as 
\begin{equation}
    \mathcal A_z^\text{adv}( z_h, v_h;\mathbf u_h):=  \sum_{e=1}^{E_h}\int_{K_e} \rho_f z_h \mathbf u_h \cdot \nabla v_h\,\mathrm d\mathbf x - \sum_{e=1}^{E_h} \int_{\partial K_e^+} \rho_f z_h^+ \mathbf u_h \cdot \mathbf n_e\,v_h\,\mathrm ds - \sum_{\gamma \in \Gamma_h^\text{in}}\int_\gamma \rho_f z_\text{in} v_h\,\mathrm ds.
\end{equation}
Here $\partial K_e^+:=\{\mathbf x\in \partial K_e : \mathbf u_h(\mathbf x)\cdot \mathbf n_e(\mathbf x) \geq 0\}$. The upwind value of ${z}_h = z^c_h + \delta z_h\in V^\text{CG}_h\oplus V^\text{DG}_h$ is given by
\begin{equation}
{z}_h^+(\mathbf x) :=
    \begin{cases}
        z_h|_{K_e}(\mathbf x) = z^c_h(\mathbf x)+\delta z_e\quad &\text{if }\mathbf x \in \partial K^+_e,\\
        z_h|_{K_{e'}}(\mathbf x) = z^c_h(\mathbf x) + \delta z_{e'}\quad &\text{if }\mathbf x \in (\partial K_e\backslash\partial K_e^+)\cap \partial K_{e'}.
    \end{cases}
\end{equation}




\subsection{Flux-corrected transport EG scheme on concentration}\label{sec:FCT}
We adopt the flux-corrected transport (FCT) scheme developed by Kuzmin, Hajduk, and Rupp~\cite{Kuzmin-EG} for an EG discretization of the linear advection equation. 
The method combines a continuous $\mathbb{Q}_1$ Galerkin approximation with a piecewise constant $\mathbb{Q}_0$ enrichment to ensure local conservation and improved stability. 
In our implementation, the FCT limiter employs the predictor--corrector formulation described in \cite[Sections~4--5]{Kuzmin-EG}.

For completeness, we recall that the EG solution admits the natural decomposition
$$z_h = \sum_{j=1}^{N_h}z^c_j \varphi_j + \sum_{e=1}^{E_h}\delta z_e \chi_e,$$
where $\{\varphi_j\}_{j=1}^{N_h}$ is the basis function of $V_h^\text{CG}$ and $\{\chi_e\}_{e=1}^{E_h}$ is the basis function of $V_h^\text{DG}$. The degrees of freedom $\{z^c_j\}_{j=1}^{N_h}$ and $\{\delta z_e\}_{e=1}^{E_h}$ define a continuous Galerkin component $z^c_h = \sum_{j=1}^{N_h}z^c_j \varphi_j$ and a piecewise-constant correction $\delta z_h=\sum_{e=1}^{E_h}\delta z_e \chi_e$. The coefficients
$z^c_j,\ j=1,\ldots,N_h$ represent the nodal values of $z^c_h$. The coefficients 
$\delta z_e$ of the discontinuous enrichment $\delta z_h$ are given by
\begin{equation}\label{eq:EG-decomposition}
\delta z_e
    = \bar{z}_e - \bar{z}^c_e,\qquad e=1,\ldots,E_h.
\end{equation}
Here $\bar{z}_e = |K_e|^{-1}\int_{K_e} z_h^{\mathrm{EG}}\,\mathrm{d}x$ and 
$\bar{z}^c_e = |K_e|^{-1}\int_{K_e} z^c_h\,\mathrm{d}x$ represent the cell averages of the EG approximation and of its CG component, respectively. 

The semi-discrete system for the coefficients of $(\bar{z}_h,z^c_h)$ can be written
in the partitioned form
\begin{equation}\label{eq:semi-egfct}
\begin{bmatrix}
    \bar M & 0\\
    0      & M_L
\end{bmatrix}
\dfrac{\mathrm d}{\mathrm dt}
\begin{bmatrix}
     \bar{z}\\[0.2em]
     z^c
\end{bmatrix}
=
\begin{bmatrix}
    \bar{a}^{\mathrm{adv}}_z(\bar{z},z^c)\\[0.2em]
    {a}^{\mathrm{adv}}_z(\bar{z},z^c) + \dot{a}_z(\bar{z},z^c)
\end{bmatrix},
\end{equation}
where $\bar{z} = (\bar{z}_e)_{e=1}^{E_h}$ is the vector of EG cell averages, 
$z^c = (z_j^c)_{j=1}^{N_h}$ is the vector of CG degrees of freedom, 
$\bar M = \mathrm{diag}(|K_e|)_{e=1}^{E_h}$ is the diagonal mass matrix of the finite volume scheme for evolving the cell averages $\bar z_e$, 
and $M_L = \mathrm{diag}(m_i)_{i=1}^{N_h}$ is the lumped mass matrix of the space $V_h^\text{CG}$, with 
positive diagonal entries $m_i = \sum_{e\in\mathcal E_i}\int_{K_e}\varphi_i\,\mathrm{d}x$. 

The right-hand side vectors 
$\bar{a}^{\mathrm{adv}}_z\in \mathbb R^{E_h}$ 
and 
${a}^{\mathrm{adv}}_z, \dot{a}_z\in \mathbb R^{N_h}$ 
are composed from
\begin{align}
    &\bar{a}^\mathrm{adv}_{z,e} := \mathcal A_z^\mathrm{upw}(z_h,\chi_e), 
    && e = 1,\ldots, E_h,\\
    &{a}^\mathrm{adv}_{z,i} := \mathcal A_z^\mathrm{adv}(z_h,\varphi_i), 
    && i = 1,\ldots, N_h
\end{align}
and
\begin{equation}
    \dot a_{z,i}:= \sum_{e\in\mathcal E_i}\int_{K_e} \varphi_i\left[
    \left(\frac{\mathrm dz^c_i}{\mathrm dt} - \frac{\mathrm dz^c_h}{\mathrm dt}\right) 
    - \left(\frac{\mathrm d\bar{z}_e}{\mathrm dt} - \frac{\mathrm d\bar{z}^c_e}{\mathrm dt}\right)
    \right]\mathrm{d}x, \qquad i = 1,\ldots, N_h.
\end{equation}
The contribution of $\dot a_{z,i}$ compensates the error due to mass lumping on
the left-hand side of the subsystem for the CG degrees of freedom \cite{Kuzmin-EG}.
The operator $\mathcal A^\text{upw}_{z}$ approximates the advective term of the subsystem for
$\bar z$  using a finite volume scheme with upwind fluxes \cite{Kuzmin-EG}.

In the predictor step of the FCT algorithm, the advective contributions are approximated by low-order operators that yield vectors
$\bar a_z^L \in \mathbb R^{E_h}$ and $a_z^L \in \mathbb R^{N_h}$ with 
\begin{align}
    &\bar a_{z,e}^L := \tilde{\mathcal A}^\text{upw}_z(z_h,\chi_e), && e = 1,\ldots,E_h,\\
    &a^L_{z,i} := \tilde{\mathcal A}^\text{adv}_{z}(z_h,\varphi_i)
        + \sum_{j\in\mathcal N^e}d^e_{ij}(z^c_j - z^c_i), && i = 1,\ldots,N_h.
\end{align}
The operator $\tilde{\mathcal A}^\text{upw}_{z}$ of the low-order subproblem for $\bar z$
corresponds to a monotone finite volume scheme in which the upwind (local Lax--Friedrichs)
fluxes use cell averages instead of the one-sided limits
of the piecewise-linear discontinuous EG approximation \cite{Kuzmin-EG,KUZMIN2025114323}.
The operator $\tilde{\mathcal A}^\text{adv}_{z}$ of the low-order subproblem for $z^c$
is constructed by removing the contribution of $\bar z$ from  $\mathcal A^\text{adv}_{z}$.
The artificial diffusion coefficients $d^e_{ij}$ depend on
$k^e_{ij} := -\int_{K_e}\varphi_i\,\mathbf u_h \cdot \nabla \varphi_j\,\mathrm{d}x$. 
The local extremum diminishing (LED) property  of the low-order CG scheme
is guaranteed for
\begin{equation}
    d^e_{ij} := 
    \begin{cases}
        \max\{-k_{ij}^e,\,0,\,-k^e_{ji}\}, & \text{if } j\in\mathcal N^e\setminus\{i\},\ i\in\mathcal N^e,\\[0.3em]
        - \sum_{k\in\mathcal N^e\setminus\{i\}}d^e_{ik}, & \text{if } j=i\in\mathcal N^e,\\[0.3em]
        0, & \text{otherwise}.
    \end{cases}
\end{equation}

The corrector step of the FCT scheme corrects the low-order predictor by adding
limited counterparts of the anti-diffusive correction terms 
$\bar b_z \in \mathbb R^{E_h}$ and $b_z \in \mathbb R^{N_h}$ consisting of
\begin{align}
    &\bar b_{z,e} := \bar a^\text{adv}_{z,e} - \bar a^L_{z,e}, && e = 1,\ldots,E_h,\\
    &b_{z,i} :=  a^\text{adv}_{z,i} -  a^L_{z,i}, && i = 1,\ldots,N_h.
\end{align}
All relevant implementation details---including the limiting procedures that ensure the LED property for forward Euler stages of an explicit SSP Runge--Kutta time-stepping method---can be found in 
\cite{Kuzmin-EG}. An extension to nonlinear scalar conservation laws was proposed in ~\cite{KUZMIN2025114323}.

To avoid spurious oscillations within the local bounds of the FCT constraints, a high-order linear or nonlinear stabilization term can be incorporated into the raw antidiffusive element contributions to be limited in  the correction stage of the FCT algorithm proposed in~\cite{Kuzmin-EG}. Example of suitable local stabilization operators can be found, e.g., in \cite[Chapter 4]{kuzmin2023}.

\subsection{IMPES framework and sequential solution strategy}
\label{sec:IMPES}

To solve the coupled flow–transport–heat system for the coefficients of the finite element approximations 
$p_h$,$\mathbf u_h$, $z_h^n$, $T_{i,h}^n$,  $T_{r,h}^n$, $ T_{s,h}^n$, we use an IMPES-type (Implicit Pressure, Explicit Concentration and Implicit Heat) scheme. To advance the numerical solution 
in time from $t^n$ to $t^{n+1}$ until a final time $t_{\mathrm{final}}$, our 
algorithm updates the quantities of interest sequentially. An implicit pressure 
update is followed by the computation of Darcy's velocity,
an explicit update of the species concentration, and implicit
updates of the temperature fields. This procedure allows us to decompose the coupled 
nonlinear system into subproblems that are solved in a decoupled yet consistent manner.
The segregated solution strategy reduces the computational cost significantly while 
maintaining numerical stability under realistic time-step constraints.

The overall workflow of the IMPES-type algorithm, as presented in Fig.~\ref{fig:flowchart}, 
illustrates the order in which individual subproblems are solved within a single time step. Repeating
the sequence of sequential updates in an iterative manner would yield the solution of
the coupled problem. However, we update each variable just once per time step. This corresponds to
using a consistent IMEX linearization / operator splitting for the global nonlinear system.

\begin{figure}[h!]
\centering
\begin{tikzpicture}[
  node distance=14mm,
  block/.style={draw, rounded corners=2mm, minimum width=34mm, minimum height=10mm, align=center},
  arrow/.style={-Latex, thick},
  dashedarrow/.style={-Latex, dashed},
  note/.style={font=\footnotesize, align=center}
]

\node[block] (p) {Darcy Flow\\($p_h,\mathbf{u}_h$)\\(Implicit)};
\node[block, below=of p] (c) {Transport\\($z_h^{n+1}$)\\(Explicit: SSP RK2)};
\node[block, below=of c] (t) {Heat\\($T_{i,h}^{n+1}, T_{r,h}^{n+1}, T_{s,h}^{n+1}$)\\(Implicit: BDF2)};


\draw[arrow] (p.south) .. controls +(0,0) and +(0,0).. node[note, left] {given $z^{n}_h$} (c.north);
\draw[arrow] (c.south) .. controls +(0,0) and +(0,0).. node[note, left] {given $T^{m}_{i,h},T^{m}_{i,h},T^{m}_{i,h}$\\
$m = n,n-1$} (t.north);
\draw[->] 
(t.east) -- ++(1cm,0) 
|- 
node[pos=0.2, right, align=left] 
{advance time\\
$t^n \rightarrow t^{n+1}$}
(c.east);



\end{tikzpicture}
\caption{Flowchart of the sequential IMPES-type procedure (update from $t^n$ to $t^{n+1}$): 
implicit Darcy update followed by an explicit transport step for the concentration and umplicit temperature updates. 
}
\label{fig:flowchart}
\end{figure}
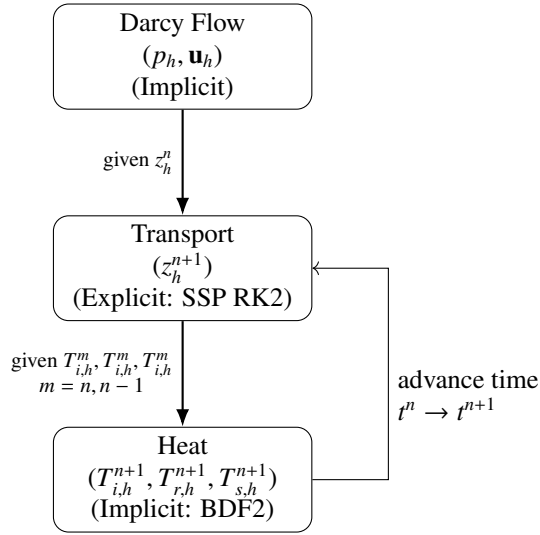

\section{Numerical Results} \label{sec:numericsim}

In this section, we present three numerical examples to demonstrate the accuracy and applicability of the proposed approach. The first example is a convergence study that validates the algorithmic coupling of the concentration equation with the three-way LTNE model; see Figure~\ref{fig:flowchart}. The second example compares the classical LTE and LTNE models for different heat transfer regimes, both with and without advection, in order to verify the asymptotic analysis presented in Section~\ref{subsec:asymptotic}. In addition, the three-way LTNE model is simulated under the same setup for comparison. The third example applies the three-way LTNE model to flow and heat transport in a fractured medium. All numerical experiments are implemented using the \texttt{deal.II} \cite{dealII94} finite element library.
\subsection{Example 1: Convergence Test}
For the first example, we test the \textit{concentration} $(z)$ and \textit{temperature} $(T_i,T_r,T_S)$ models with respect to prescribed smooth solutions in the two-dimensional domain $\Omega = [0,1]^2$ with a final computation time $t_\mathrm{final} = 0.1$, defined as follows:
\begin{equation*}
\begin{split}
    z(t,x,y) &= \sin(t)\,\sin\!\big(\pi (x+y)\big)\,\cos\!\big(\pi(2x - y)\big),\\
    T_i(t,x,y) &= \cos(t)\,\cos\!\big(2\pi (x+y)\big)\,\sin\!\big(\pi(2x - y)\big),\\
    T_r(t,x,y) &= \sin(t)\,\sin\!\big(2\pi (x+y)\big)\,\cos\!\big(\pi(2x - y)\big),\\
    T_s(t,x,y) &= t\,\sin(x+y).
\end{split}
\end{equation*}
A prescribed exact divergence-free velocity field is given by
\begin{equation*}
    \mathbf{u}(x,y) =
    \left(
        \sin(\pi x)\cosh\!\big(\pi(y-\tfrac{1}{2})\big),\;
        -\cos(\pi x)\sinh\!\big(\pi(y-\tfrac{1}{2})\big)
    \right).
\end{equation*}
For simplicity, the physical parameters in the energy conservation laws are assumed to be constant
\begin{equation*}
\begin{split}
    \rho_\beta = c_\beta = 1, \quad 
    \lambda_{\beta,\mathrm{eff}} = 0.5, 
    &\qquad \beta \in \{i,r,s\},\\
    h_{\beta_1\beta_2} = 0.1, 
    &\qquad \beta_1 \neq \beta_2 \in \{i,r,s\}.
\end{split}
\end{equation*}
The sources are prescribed to satisfies the equations \eqref{eq:z},\eqref{eq:Ti}-\eqref{eq:Ts}. The Dirichlet value of $z,T_i,T_r$ are prescribed on in flow boundary ($\Gamma_\text{in} = \{(x,y)\in\partial \Omega: x=1\}$), and exact Neumann values for $T_i,T_r,T_s$ are prescribed on the whole boundary $\partial \Omega$.

For each time step, we first solve concentration $z$ using EG-FCT method using SSP-RK2, then solve all of the temperature $(T_i,T_r, T_s)$ use EG-$\mathbb Q_1$ with BDF2 time discretization, and
the problem are solved on a sequence of four uniformly refined meshes ranging from $h = 2^{-3}$ to $h = 2^{-6}$ with an initial time step of $\Delta t = 0.001$, which is successively halved at each refinement cycle.

The convergence results are presented in Table~\ref{tab:convergency}, and the 
errors are evaluated using the following formulas, demonstrated here for the variable $z$:
\begin{equation}
    \max_{0 \leq t^n \leq t_{\text{final}}} 
    \| z_h^n - z^n \|_{L^2(\Omega)}, 
    \qquad
    \max_{0 \leq t^n \leq t_{\text{final}}} 
    \| z_h^n - z^n \|_{H^1(\mathcal{T}_h)}.
\end{equation}
From Table~\ref{tab:convergency}, we observe that the convergence rate in the $L^2$-norm increases with mesh refinement, approaching nearly second-order accuracy, while the convergence rate in the $H^1$-norm attains first-order accuracy. These results demonstrate the expected optimal convergence behavior of the proposed numerical scheme.

\begin{table}[H]
\centering
\begin{tabular}{c c c c c c}
\hline
\multirow{2}{*}{Variable} & \multirow{2}{*}{Refinement}
& \multicolumn{2}{c}{$\max_n \|\cdot\|_{L^2}$}
& \multicolumn{2}{c}{$\max_n \|\cdot\|_{H^1}$} \\ 
\cline{3-6}
 &  & Error & Rate & Error & Rate \\ 
\hline
\multirow{4}{*}{$T_i$}
& 1 & 4.8666e$-$01 & --   & 4.1301e$+$00 & --   \\
& 2 & 1.2136e$-$01 & 2.00 & 2.2812e$+$00 & 0.86 \\
& 3 & 3.0415e$-$02 & 2.00 & 1.1598e$+$00 & 0.98 \\
& 4 & 7.6041e$-$03 & 2.00 & 5.7955e$-$01 & 1.00 \\
\hline
\multirow{4}{*}{$T_r$}
& 1 & 1.5885e$-$02 & --   & 2.0526e$-$01 & --   \\
& 2 & 3.9710e$-$03 & 2.00 & 1.1285e$-$01 & 0.86 \\
& 3 & 9.9281e$-$04 & 2.00 & 5.7571e$-$02 & 0.97 \\
& 4 & 2.4820e$-$04 & 2.00 & 2.8900e$-$02 & 0.99 \\
\hline
\multirow{4}{*}{$T_s$}
& 1 & 1.3053e$-$03 & --   & 4.2931e$-$03 & --   \\
& 2 & 3.2632e$-$04 & 2.00 & 2.2297e$-$03 & 0.95 \\
& 3 & 8.1581e$-$05 & 2.00 & 1.0687e$-$03 & 1.06 \\
& 4 & 2.0395e$-$05 & 2.00 & 5.1993e$-$04 & 1.04 \\
\hline
\multirow{4}{*}{$z$}
& 1 & 1.5885e$-$02 & --   & 2.0524e$-$01 & --   \\
& 2 & 3.9702e$-$03 & 1.98 & 1.1245e$-$01 & 0.87 \\
& 3 & 9.9241e$-$04 & 1.99 & 5.7571e$-$02 & 0.98 \\
& 4 & 2.4820e$-$04 & 2.00 & 2.8700e$-$02 & 0.99 \\
\hline
\end{tabular}

\caption{Convergence results for $\max_n \|\cdot\|_{L^2}$ and $\max_n \|\cdot\|_{H^1}$ norms for $z$, $T_i$, $T_r$, and $T_s$.}
\label{tab:convergency}
\end{table}

\begin{remark}
We emphasize that the EG-FCT limiting procedure is applied only to the advection fluxes and does not modify the source terms in the concentration equation.
In this convergence test, the prescribed solution is smooth; consequently, the limiter remains almost inactive or only weakly active (with the flux limiter almost to be 1).
As a result, the method retains its underlying high-order accuracy, and the expected optimal convergence rates are observed.
\end{remark}

\subsection{Example 2: Comparison between LTE and LTNE models}
In the second example, we want to discuss where the LTE assumption is insufficient and the LTNE model therefore becomes necessary. We consider a computational domain 
$\Omega = [0,4]\times[0,1]$ with a final simulation time of $t_\text{final} = 1.0$. 
 We will test LTNE problem with and without advection, combined with different ranges of the heat transfer coefficient $h_{fs}$, while the volume-specific interfacial area is fixed at $a_{fs} = 1$. The following two cases are considered:
\begin{itemize}
    \item Case 1: LTNE with fast heat exchange ($h_{fs} = 10$)  
    \begin{itemize}
        \item Without and with advection ($\mathbf{u} = (0.5,0)$).  
        \item Compared with LTE with Dirichlet condition leans toward the fluid phase \eqref{eq:example1-bc-fast-heat chenge}.  
    \end{itemize}
    \item Case 2: LTNE with slow heat exchange ($h_{fs} = 0.001$)  
    \begin{itemize}
        \item Without and with advection ($\mathbf{u} = (0.5,0)$).  
        \item Compared LTE with Dirichlet condition leans toward the phase average \eqref{eq:example1-bc-slow-heat chenge}  
    \end{itemize}
\end{itemize}
Details of the boundary and initial conditions are shown in Figure~\ref{fig:domain-setup for LTNE example1}.


For simplicity, we neglect the source terms by setting $Q_f = Q_s = 0$, and assume a porosity of $\phi = 0.5$. 
Both phases are taken to have identical densities, heat capacities, 
and conductivities 
$$\rho_f = \rho_s = 1,\quad  c_f = c_s = 1,\quad \lambda_f = \lambda_s = 1.$$
We also set the same initial conditions for solid and fluid phases as 
\begin{equation}
    T_f = T_s = 100 \quad \text{on }\Omega \times\{0\},
\end{equation}
and impose the following boundary conditions 
\begin{equation}
\begin{split}
    T_f  = 20 \quad &\text{on }\Gamma_\mathrm D \times (0,t_\text{final}],\\
   \lambda^\text{eff}_f \nabla T_f \cdot \mathbf n = 0\quad &\text{on }\partial\Omega \backslash\Gamma_D \times(0,t_\text{final}],\\
   \lambda_s^\text{eff} \nabla T_s\cdot \mathbf n = 0 \quad &\text{on }\partial\Omega\times(0,t_\text{final}],
\end{split}
\end{equation}
where $\Gamma_D = \{(x,y)\in\partial\Omega: x=0\}$. Although the problem is posed in a two-dimensional domain, the solutions are essentially one-dimensional. 


\begin{figure}
\centering
\begin{tikzpicture}[scale=1.4, every node/.style={font=\footnotesize}]
  \def\L{4}  
  \def\H{1}  

  \tikzset{
    dirEdge/.style={line width=1.3pt,blue},
    neuEdge/.style={line width=0.8pt,black},
    normalvec/.style={-Latex, line width=0.6pt},
    meas/.style={Latex-Latex, line width=0.5pt},
    neuMark/.style={decoration={markings, mark=at position 0.5 with {\arrow{Stealth}}},
                    postaction=decorate}
  }

  \draw[neuEdge] (0,0) rectangle (\L,\H);

  \draw[dirEdge] (0,0) -- (0,\H)
    node[pos=0.9, left=4pt, black] {Dirichlet Boundary $\Gamma_{\mathrm D}$  };
    \draw[dirEdge] (0,0) -- (0,\H)
    node[pos=0.65, left=8pt, black] {$T_f = 20$};
    \draw[dirEdge] (0,0) -- (0,\H)
    node[pos=0.3, left=8pt, black] {Case 1: $T = 20$};
    \draw[dirEdge] (0,0) -- (0,\H)
    node[pos=0.1, left=8pt, black] {Case 2: $T = 60$};

  \draw[neuEdge] (\L,0) -- (\L,\H)
    node[pos=0.9, right=2pt, black] {Neumann boundary $\partial\Omega \backslash\Gamma_\mathrm D$};
     \draw[neuEdge] (\L,0) -- (\L,\H)
    node[pos=0.65, right=4pt, black] {$\lambda^\text{eff}_f\nabla T_f \cdot \mathbf n = \lambda_f^\text{eff}\cdot \mathbf n = 0$};
    \draw[neuEdge] (\L,0) -- (\L,\H)
    node[pos=0.23, right=4pt, black] {Case 1\& 2: $\lambda_\text{eff}\nabla T \cdot \mathbf n = 0$};

  \node at ($(0,-0.5)!0.5!(\L,\H)$) {$\Omega=[0,4]\times[0,1]$};
  \node at ($(0,0)!0.5!(\L,\H)$) {Initial: $T_f=T_s = T= 100$};

\end{tikzpicture}
\caption{Domain and boundary/initial conditions for LTNE models ($T_f$ and $T_s$), and LTE ($T$) for Case 1 and Case 2.}
\label{fig:domain-setup for LTNE example1}
\end{figure}
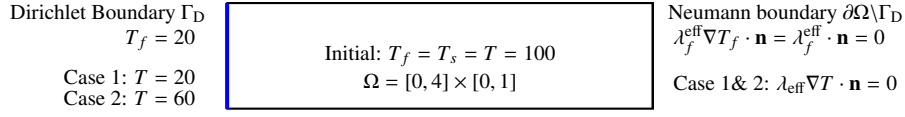


\subsubsection*{LTNE: fast heat exchange $h_{fs} = 10$ between $T_f$ and $T_s$:}
We first consider the case $h_{fs} = 10$, which represents scenarios where the heat 
exchange between $T_f$ and $T_s$ is fast (see Figure~\ref{fig:xxx-without advection}-\ref{fig:xxx-withadvection}). 
The LTNE problem is solved both with advection, $\mathbf{u} = (0.5,0)$, and without advection, 
using the numerical scheme described in Section~\ref{sec:numerics}, 
on a mesh of size $h = 2^{-5}$ with a timestep $\Delta t = 0.005$. 
For comparison, we also solve the LTE model with the following initial and boundary conditions:
\begin{equation}\label{eq:example1-bc-fast-heat chenge}
    T = 100 \quad \text{on } \Omega \times \{0\}, \qquad 
    \begin{cases}
        T = 20 & \text{on } \Gamma_D \times (0,t_\text{final}], \\[6pt]
        \lambda_\text{eff}\nabla T \cdot \mathbf{n} = 0 & \text{on } \partial\Omega \times (0,t_\text{final}].
    \end{cases}
\end{equation}
Note that the Dirichlet boundary condition for the equilibrium temperature $T$ 
is chosen so that it leans toward the fluid temperature $T_f$ on $\Gamma_D$.

We present the solution values along the middle horizontal line 
$\{(x,y) \in \Omega : y = 0.5\}$. The numerical solutions are shown in 
Figure~\ref{fig:xxx-without advection} (without advection) and 
Figure~\ref{fig:xxx-withadvection} (with advection), where the mixed temperature 
$T_\text{mix} = \tfrac{1}{2}(T_f + T_s)$, computed from the LTNE model as defined in 
\eqref{eq:Tmix-def}, is also included. In both cases, the two temperatures $T_f$ and $T_s$ 
rapidly approach each other and reach the same value within a short time. As a result, we 
also observe $T_\text{mix} \approx T$, which is consistent with our analysis in 
Section~\ref{subsec:asymptotic}.

 \begin{figure}
\centering
    \includegraphics[width=0.8\linewidth]{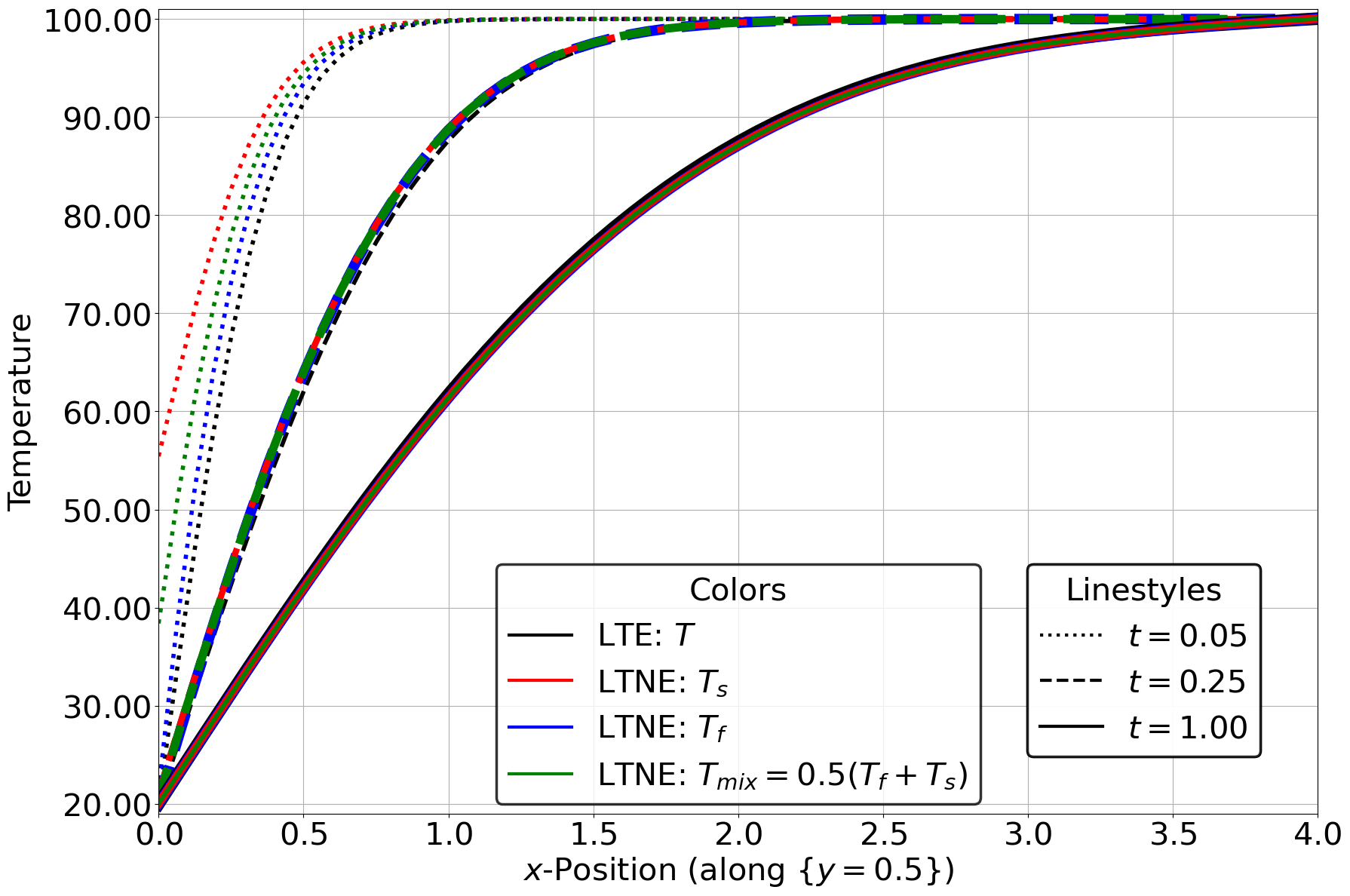}

    \caption{Example 1: LTNE: fast heat exchange $h_{fs} = 10$, without advection. The solution values of $T_f, T_s$, $T_\text{mix}= 0.5(T_f+T_s)$, and $T$ from LTE with initial/boundary condition \eqref{eq:example1-bc-fast-heat chenge}. }
    \label{fig:xxx-without advection}
\end{figure}

 \begin{figure}
 \centering
\includegraphics[width=0.8\linewidth]{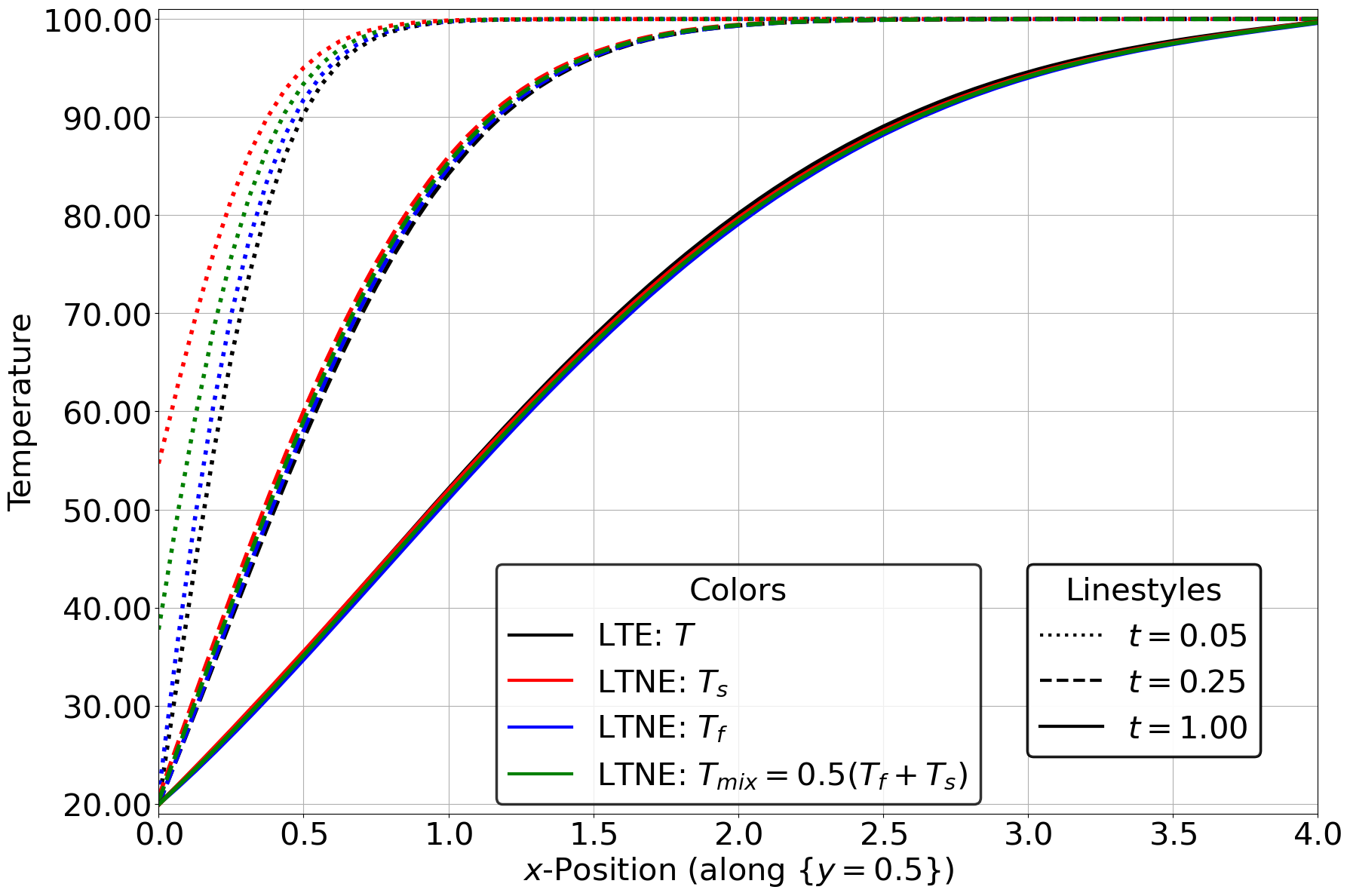}
    \caption{Example 1: LTNE: fast heat exchange $h_{fs} = 10$, with advection $\mathbf u = (0.5,0)$. The solution values of $T_f, T_s$, $T_\text{mix}= 0.5(T_f+T_s)$, and $T$ from LTE with initial/boundary condition \eqref{eq:example1-bc-fast-heat chenge}. }
    \label{fig:xxx-withadvection}
\end{figure}

\subsubsection*{LTNE: Slow heat exchange $h_{fs} = 0.001$ between $T_f$ and $T_s$:}
We next set $h_{fs} = 0.001$, representing scenarios where the heat exchange 
between $T_f$ and $T_s$ is extremely slow. The LTNE problem is solved both with advection, 
$\mathbf{u} = (0.5,0)$, and without advection, using the same discretization setup as in 
the previous case. For comparison, the LTE model is supplemented with different boundary 
conditions:
\begin{equation}\label{eq:example1-bc-slow-heat chenge}
    T = 100 \quad \text{on } \Omega \times \{0\}, \qquad 
    \begin{cases}
        T = 60 & \text{on } \Gamma_D \times (0,t_\text{final}], \\[6pt]
        \lambda_\text{eff}\nabla T \cdot \mathbf{n} = 0 & \text{on } \partial\Omega \times (0,t_\text{final}].
    \end{cases}
\end{equation}
Here, the Dirichlet boundary value for the equilibrium temperature $T$ is chosen so that 
it leans toward the equilibrium temperature between the fluid and solid, 
$\tfrac{1}{2}(T_f+T_s) = 60$, on $\Gamma_D$.

Figure~\ref{fig:slow-withoutadvection} shows the solution without advection, where we 
observe a significant difference between $T_f$ and $T_s$. Nevertheless, we still find 
$T_\text{mix} \approx T$. On the other hand, for the case with advection, shown in 
Figure~\ref{fig:slow-withadvection}, we observe not only a pronounced difference between 
$T_f$ and $T_s$, but also that $T$ no longer matches $T_\text{mix}$. This demonstrates 
that the LTE model, which assumes $T_f \approx T_s$ and approximates their equilibrium by 
a single temperature $T$, is no longer valid. In fact, it even fails to accurately 
approximate the true equilibrium temperature $T_\text{mix}$ between the two phases.

\begin{figure}
\centering
        \includegraphics[width=0.8\linewidth]{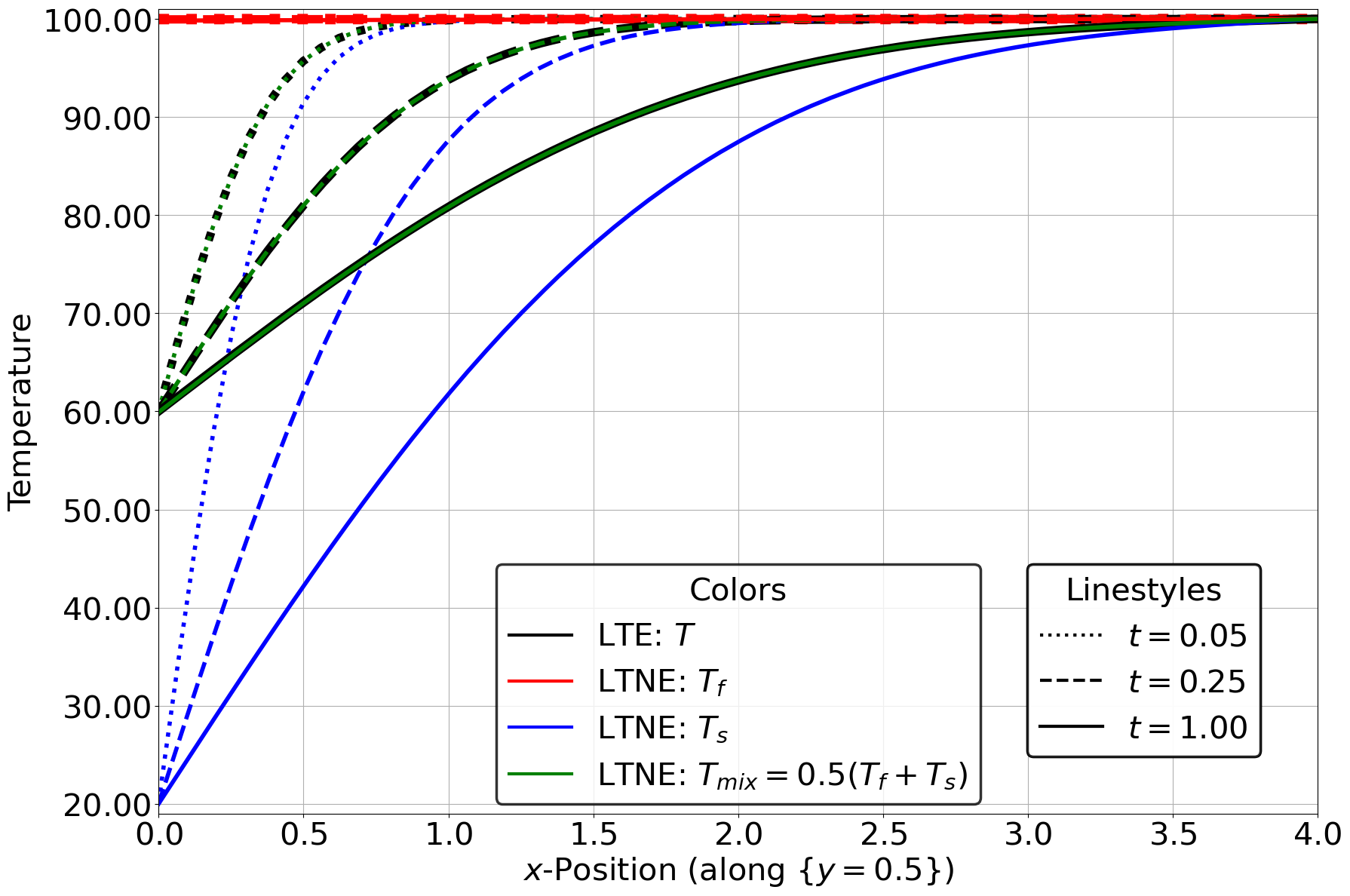}

    \caption{Example 1: LTNE: slow heat exchange $h_{fs} = 0.001$, without advection. The solution values of $T_f, T_s$, $T_\text{mix}= 0.5(T_f+T_s)$, and $T$ from LTE with initial/boundary condition \eqref{eq:example1-bc-slow-heat chenge}. }
    \label{fig:slow-withoutadvection}
\end{figure}

\begin{figure}
\centering
    \includegraphics[width=0.8\linewidth]{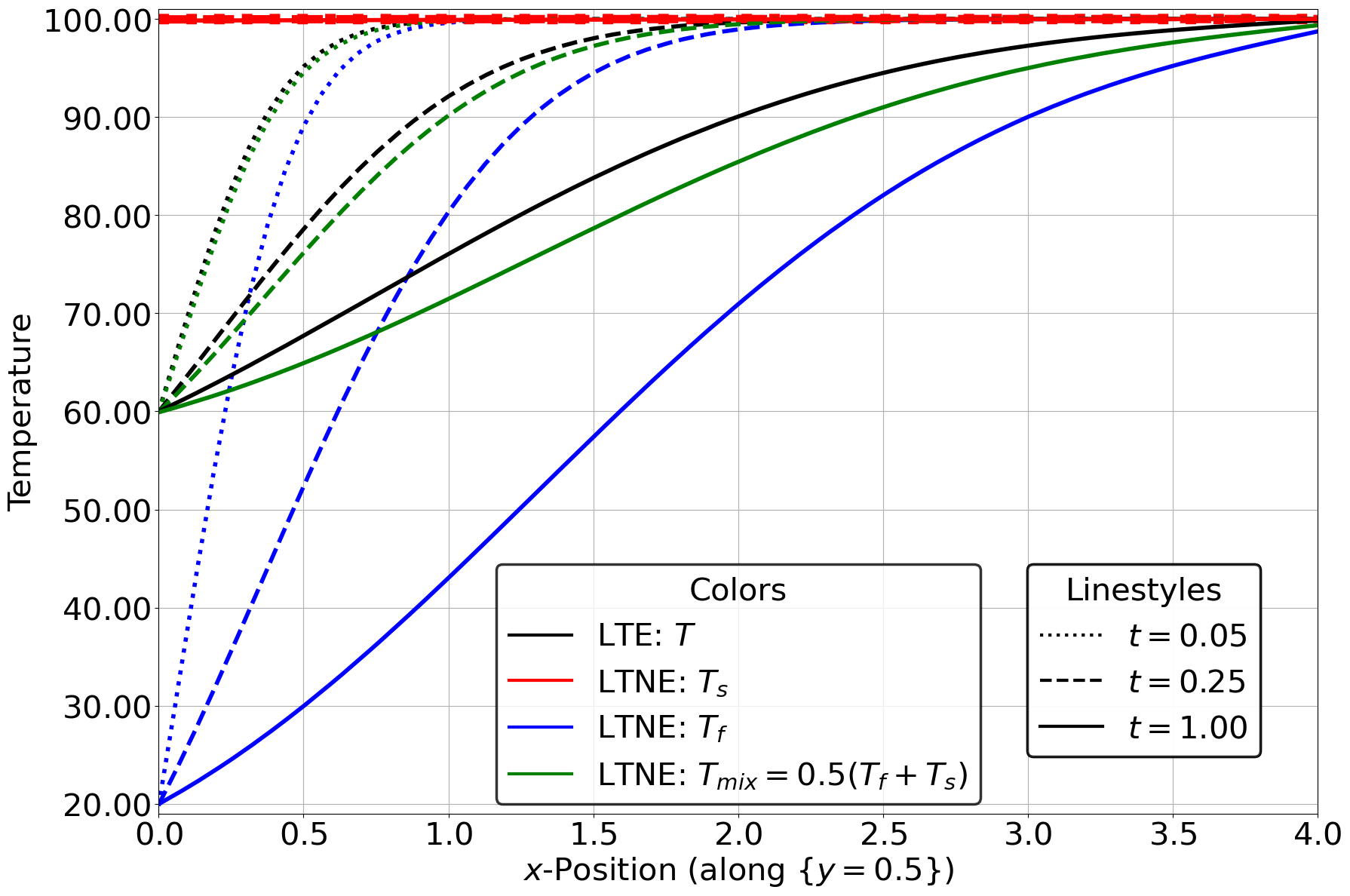}

    \caption{Example 1: LTNE: slow heat exchange $h_{fs} = 0.001$, wit advection $\mathbf u = (0.5,0)$. The solution values of $T_f, T_s$, $T_\text{mix}= 0.5(T_f+T_s)$, and $T$ from LTE with initial/boundary condition \eqref{eq:example1-bc-slow-heat chenge}. }
    \label{fig:slow-withadvection}
\end{figure}

\subsubsection*{Conclusion:} 
From the previous two setups, the temperatures $T_f$ and $T_s$ quickly approach each other 
and reach the same value within a very short time whenever $h_{fs}$ is large, regardless of 
whether advection is included or not. However, when $h_{fs}$ is small, $T_f$ and $T_s$ 
no longer reach local equilibrium. Moreover, advection amplifies the discrepancy between 
the mixed temperature $T_\text{mix} = \tfrac{1}{2}(T_f + T_s)$ and the reference temperature $T$.

We conclude the first example by plotting 
\[
\theta := T_s - T_f
\]
for various values of $h_{fs} = 1, 3, 7$, with and without advection, in 
Figure~\ref{fig:Ex1-theta}. This allows us to assess the validity of 
the LTE assumption ($T_f \approx T_s \approx T$). The results indicate that when the heat 
exchange coefficient $h_{fs}$ is not sufficiently large, the LTE assumption breaks down. 
Furthermore, the inclusion of advection exacerbates the violation of LTE, thereby 
showing the necessity of the LTNE model to describe the evolution of $T_f$, $T_s$, 
and the equilibrium temperature $T_\text{mix}$ between the two phases.

\begin{figure}
    \centering
    \includegraphics[width=1.0\linewidth]{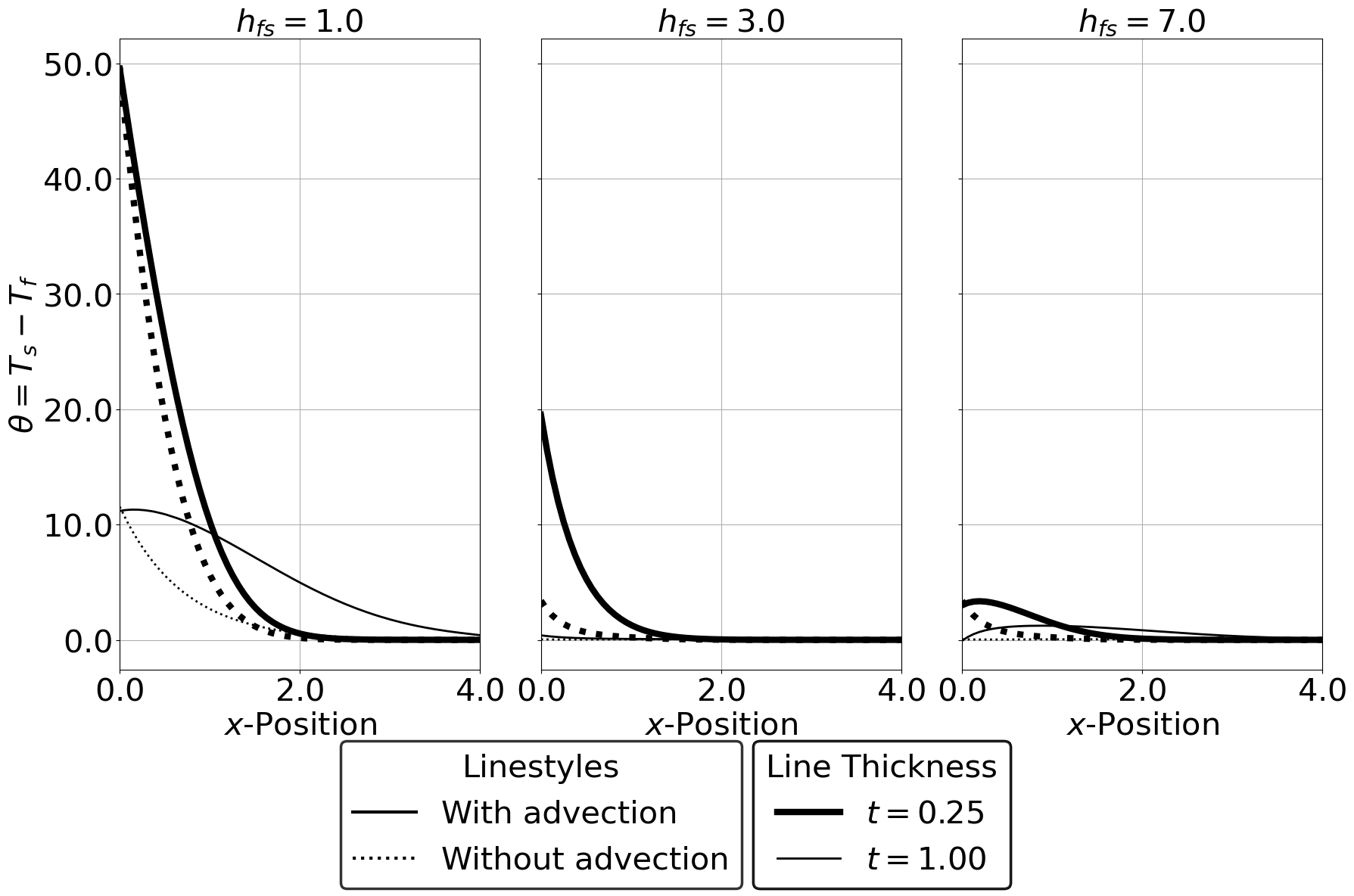}
    \caption{Example 1: Influence of heat transfer coefficient ($h_{fs}$), advection, and time on the deviation between solid and fluid temperatures $\theta = T_s - T_f$. Lower $h_{fs}$ increases the deviation, and advection further amplifies the difference between the phases and therefore, break the LTE assumption ($T_s \approx T_f \approx T$).}
    \label{fig:Ex1-theta}
\end{figure}

{
\subsection{Example 2.1: three-way LTNE Models}
In this example, we consider a three-way LTNE model with advection velocity
\[
\mathbf{u} = (0.25, 0)
\]
defined on the same domain $\Omega = [0,4] \times [0,1]$, using the same final simulation time $t_\text{final} = 8.0$. For the concentration $z$, the following inflow boundary condition and initial condition are prescribed:
\begin{equation*}
    z = 0.99 \quad \text{on } \Gamma_{\text{in}} \times (0,t_\text{final}], 
    \qquad 
    z = 0.01 \quad \text{on } \Omega \times \{0\}.
\end{equation*}

The additional parameters in the three-way LTNE model are set to $h_{\text{ir}} = a_{\text{ir}} = 1.0$, while all remaining parameters are kept the same as in the previous example. We note that a small heat exchange rate $h_{fs} = 0.001$ is chosen so that the discussion mainly focuses on the solution behavior of $T_i$ and $T_r$. Under this setting, the solid temperature $T_s$ remains nearly constant, as observed in Example~1.

\begin{figure}[H]
    \centering
    \begin{subfigure}[t]{0.8\textwidth}
        \centering
       \includegraphics[width=0.9\linewidth]{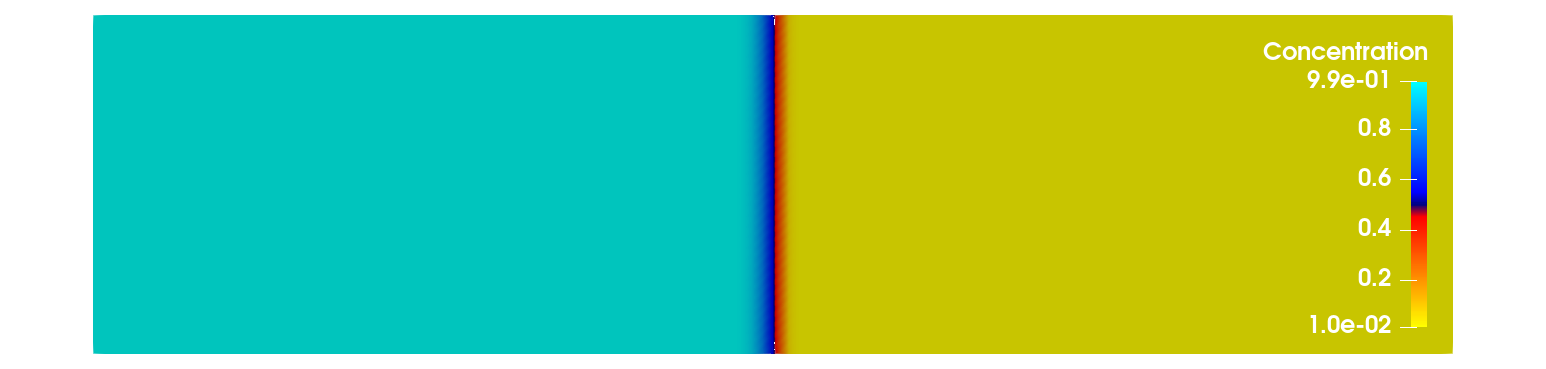}
        \caption{$z$}
        \label{example2.1-z}
    \end{subfigure}
    \begin{subfigure}[t]{0.8\textwidth}
        \centering
       \includegraphics[width=0.9\linewidth]{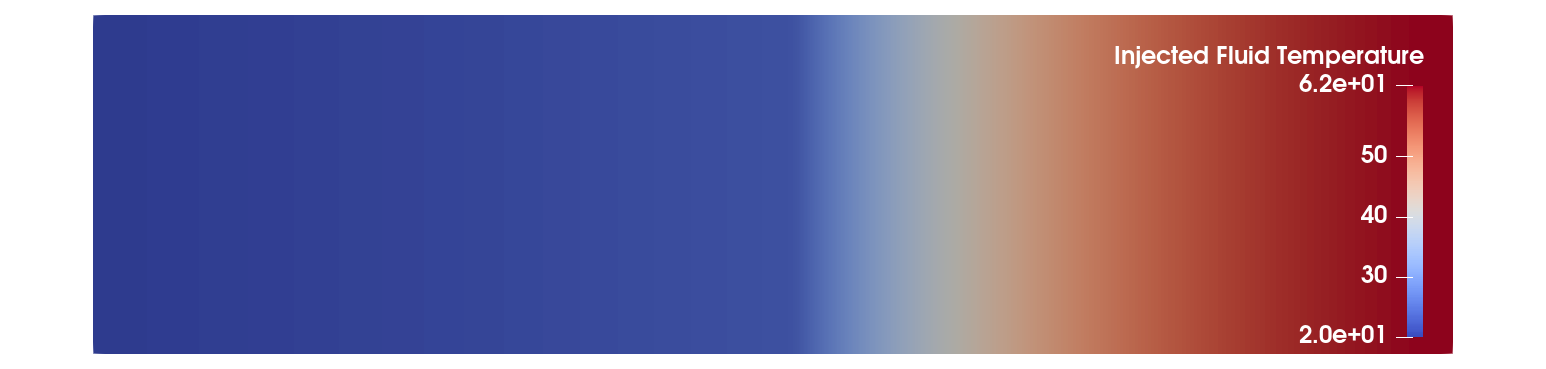}
        \caption{$T_i$}
        \label{example2.1-Ti}
    \end{subfigure}
    \begin{subfigure}[t]{0.8\textwidth}
        \centering
       \includegraphics[width=0.9\linewidth]{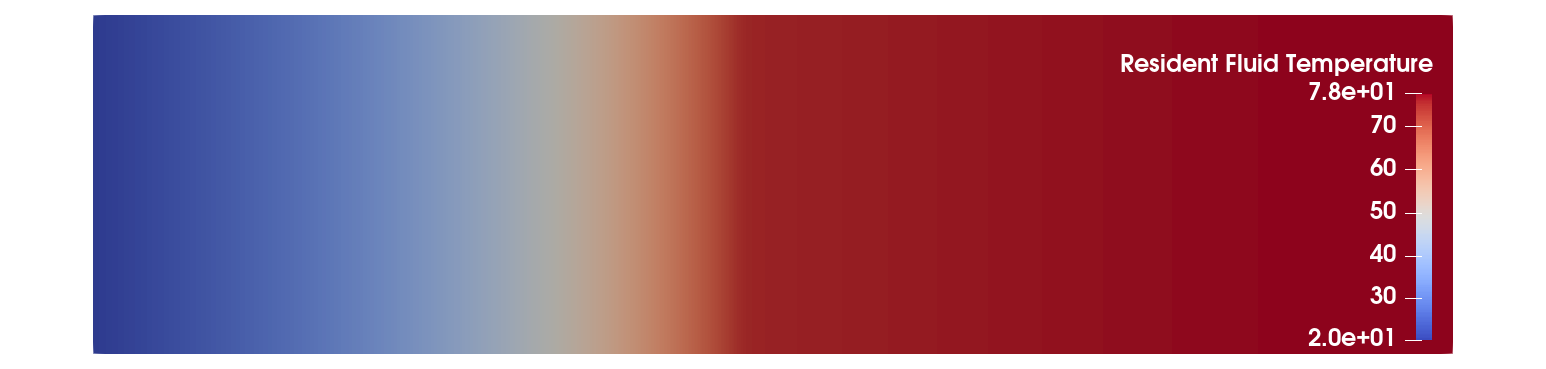}
        \caption{$T_r$}
        \label{example2.1-Tr}
    \end{subfigure}
    \begin{subfigure}[t]{0.8\textwidth}
        \centering
       \includegraphics[width=0.9\linewidth]{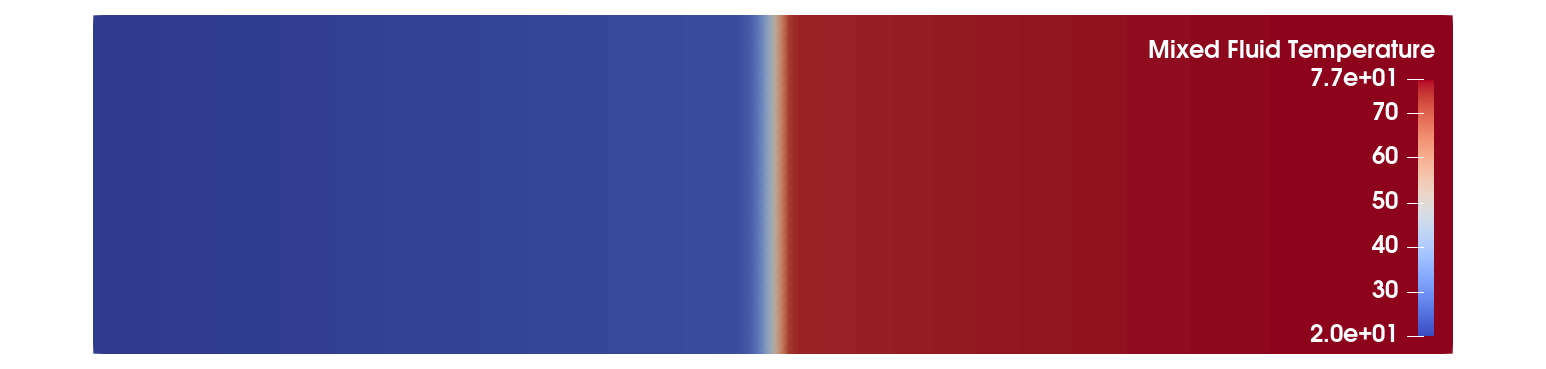}
        \caption{$T_f^{\text{mix}}$}
        \label{example2.1-Tmix}
    \end{subfigure}
    \caption{Example 2.1: The solution profiles for injected fluid temperature $T_i$, resident fluid temperature $T_r$, and computed fluid $T_f^\text{mix}$ from three-way LTNE \eqref{eq:enthalpy for fluid} and LTNE}
    \label{fig:example2.1-solution-profiles}
\end{figure}

Figures \ref{example2.1-Ti}--\ref{example2.1-Tmix} show the concentration $z$ profiles, the temperature profiles of $T_i$ and $T_r$, and the mixed fluid temperature $T_f^{\text{mix}}$ defined by \eqref{eq:enthalpy for fluid} 
\[
T_f^{\text{mix}} = z T_i + (1 - z) T_r
\]
at time $t = 8.0$. Because the setup of this example is essentially equivalent to a one-dimensional problem, the solution exhibits predominantly one-dimensional behavior.

For a more detailed examination, we plot the solution values along the middle line
$\{(x,y)\in\Omega : y = 0.5\}$
in Figures \ref{example2.temperature-overline} and \ref{example2.z-overline}, showing the temperature fields $T_i$, $T_r$, $T_s$, and $T_f^{\text{mix}}$, and the concentration $z$, respectively. As time progresses, the resident fluid temperature $T_r$ decreases because heat is extracted by the injected fluid. Consequently, the injected fluid temperature $T_i$ increases over time. Meanwhile, the mixed fluid temperature $T_f^\text{mix}$ also decreases with time. However, in contrast to the two-way LTNE model shown in Figure \ref{fig:slow-withadvection} in Example~2, the sharp concentration gradient leads to a noticeable temperature jump in this case.

\begin{figure}[H]
    \centering
    \centering
    \begin{subfigure}[t]{0.8\textwidth}
        \centering
       \includegraphics[width=0.8\linewidth]{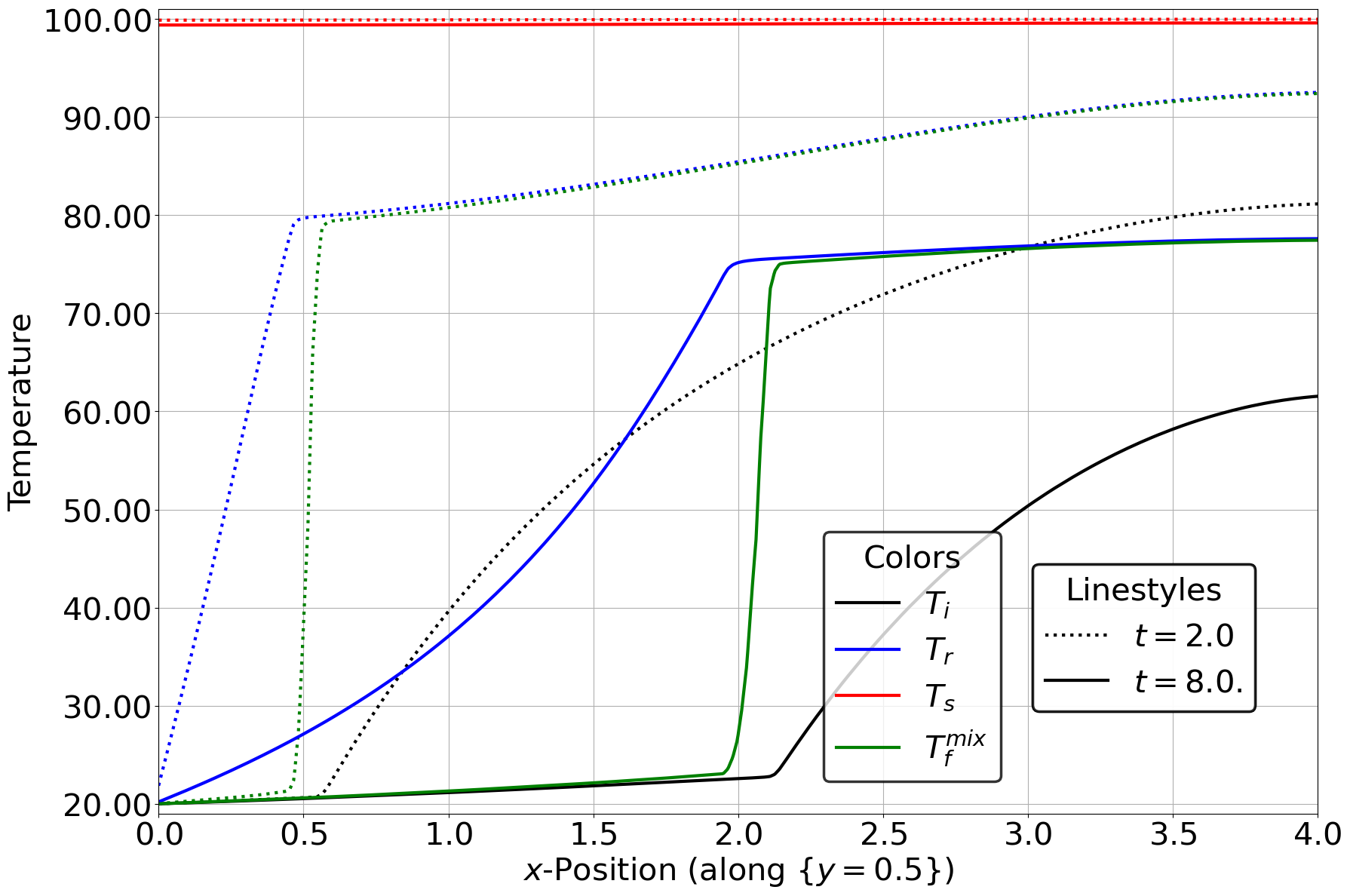}
        \caption{$T_i,T_r,T_s, T_f^\text{mix}$}
        \label{example2.temperature-overline}
    \end{subfigure}
        \begin{subfigure}[t]{0.8\textwidth}
        \centering
       \includegraphics[width=0.8\linewidth]{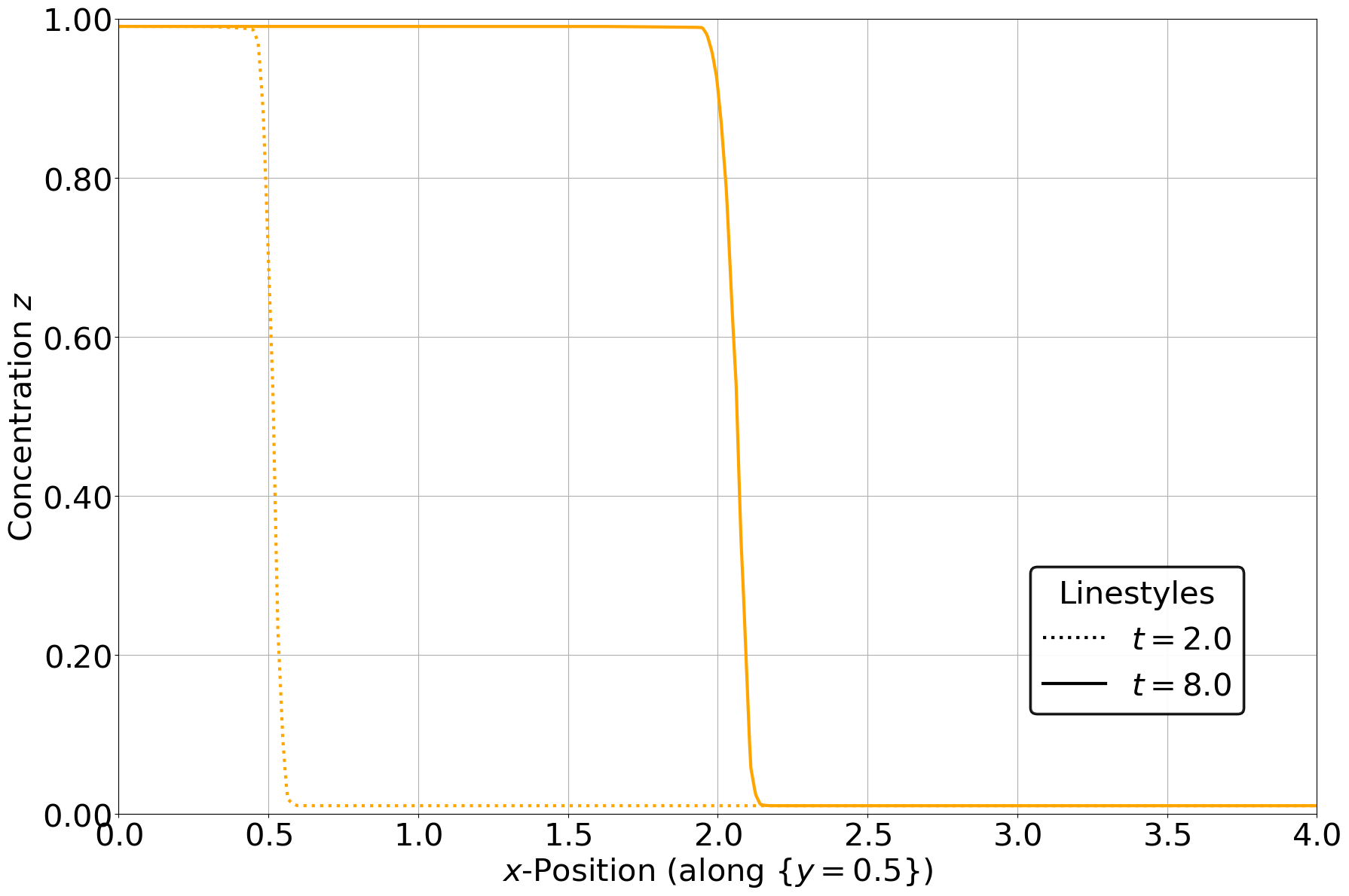}
        \caption{$z$}
        \label{example2.z-overline}
    \end{subfigure}
    \caption{Example 2.1: solution values over the middle horizontal line.}
    \label{fig:placeholder-1}
\end{figure}

}
\clearpage
\newpage



\clearpage
\newpage

\subsection{Example 3: Heterogeneous permeability with fracture: LTNE for injected and resident fluids}  
In the computational domain $\Omega = [0,4]\times[0,1]$, we consider a heterogeneous setup with a network of 34 high-permeability fracture channels (see Figure~\ref{fig:fracture-setup}). Each fracture is represented by a line segment of uniform thickness corresponding to an effective fracture width of $0.1$. The fractures are assigned a permeability of $K=1$, while the surrounding background matrix has a much lower permeability of $K=0.01$. 

We set source for solid phase to be  $Q_s\equiv 1$, and for incoming and resident fluid, we set $Q_i = Q_r \equiv  0$. The inflow boundary conditions for the energy system are prescribed as follows: 
\begin{equation}
\begin{split}
    T_i = T_r = 20 \quad &\text{on }\Gamma_{\mathrm D_1}\times(0,t_\text{final}], \\ 
    \lambda^\text{eff}_f \nabla T_i \cdot \mathbf n 
    = \lambda^\text{eff}_f \nabla T_r \cdot \mathbf n = 0 \quad 
    &\text{on } (\Gamma_{\mathrm D_2}\cup\Gamma_\mathrm N)\times(0,t_\text{final}],\\
    \lambda_s^\text{eff} \nabla T_s \cdot \mathbf n = 0 \quad 
    &\text{on }\partial\Omega\times(0,t_\text{final}].
\end{split}
\end{equation}
The boundary conditions for the concentration of incoming fluid and the pressure are given as
\begin{equation}
\begin{split}
    z = 0.99 \quad &\text{on }\Gamma_{\mathrm D_1}\times(0,t_\text{final}],\\
    p = 1 \quad &\text{on }\Gamma_{\mathrm D_1},\\
    p = 0 \quad &\text{on }\Gamma_{\mathrm D_2},\\
    \mathbf u\cdot\mathbf n = 0 \quad &\text{on }\Gamma_\mathrm N.
\end{split}
\end{equation}
The corresponding initial state for temperature and concentration are given by
\begin{equation}
    T_i = 20,\quad T_r = 100,\quad T_s = 100,\quad 
    z \equiv 0.01 
    \qquad \text{on }\Omega.
\end{equation}
The physical parameters used in the simulations are summarized in Table~\ref{tab:example2-parameters}. We use the same spatial discretization with $h = 2^{-6}$ and temporal discretization $\Delta t= 0.005$.

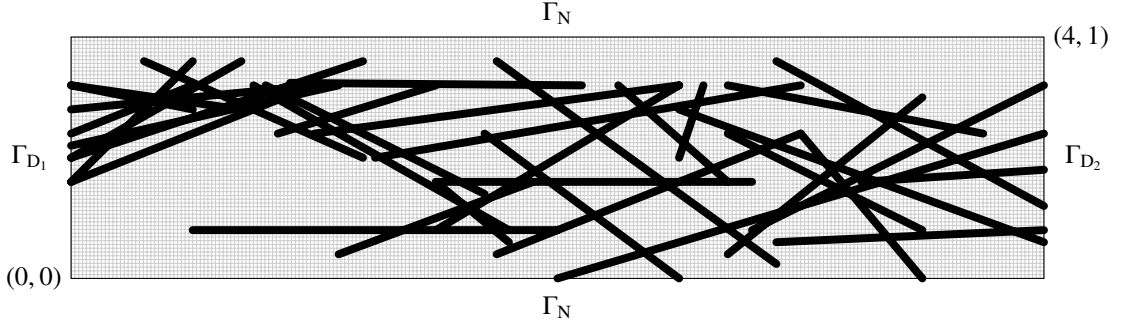
\begin{figure}[]
    \centering
    
    \begin{tikzpicture}[x=4cm,y=4cm,line cap=round,scale = 0.8]
  \def\chanlw{3.pt} 
  \def\legx{2.6}     
  \def\legy{0.88}    

  


   \draw[step=1/64, gray!40, very thin] (0,0) grid (4,1);
   \draw[] (0,0) rectangle (4,1);

  \node[left] at (-0.05,0.5) {$\Gamma_{{\mathrm D}_1}$};
  \node[right] at (4.05,0.5) {$\Gamma_{{\mathrm D}_2}$};
  \node[] at (2, 1.1) {$\Gamma_\mathrm{ N}$};
  \node[] at (2, -0.12) {$\Gamma_\mathrm{ N} $};
  \node[left] at (0,0) {$(0,0)$};
  \node[right] at (4,1) {$(4,1)$};
  
  \draw[line width=\chanlw]
    (1.5,0.2) -- (2.5,0.8)
(2.25,0.8) -- (2.7,0.4)
(1.5,0.4) -- (1.8,0.15)
(2.5,0.5) -- (2.6,0.8)
(0.0,0.8) -- (0.7,0.7)
(0.0,0.4) -- (1.0,0.8)
(0.75,0.8) -- (1.8,0.2)
(0.85,0.6) -- (1.5,0.8)
(1.75,0.9) -- (2.9,0.06)
(2.7,0.1) -- (3.5,0.75)
(2.9,0.9) -- (4.0,0.3)
(2.7,0.8) -- (3.75,0.6)
(3.25,0.4) -- (4.0,0.45)
(2.8,0.2) -- (4.0,0.8)
(0.0,0.5) -- (1.2,0.9)
(0.5,0.2) -- (2.0,0.2)
(1.5,0.4) -- (2.8,0.4)
(1.25,0.5) -- (3.0,0.8)
(1.0,0.6) -- (2.5,0.8)
(0.0,0.8) -- (0.5,0.7)
(0.0,0.6) -- (0.5,0.8)
(0.0,0.4) -- (0.5,0.9)
(0.3,0.9) -- (1.2,0.5)
(3.0,0.6) -- (3.5,0.0)
(1.7,0.6) -- (2.5,0.0)
(2.0,0.0) -- (4.0,0.6)
(2.7,0.6) -- (3.5,0.2)
(2.9,0.15) -- (4.0,0.2)
(1.75,0.1) -- (3.0,0.6)
(0.0,0.55) -- (1.1,0.8)
(0.9,0.81) -- (2.1,0.8)
(2.5,0.7) -- (4.0,0.15)
(1.1,0.1) -- (1.9,0.4)
(0.0,0.5) -- (0.7,0.9)
(0.0,0.7) -- (1.0,0.8)
(0.8,0.8) -- (1.7,0.35);
\end{tikzpicture}
    \caption{Example 2. 34 high-permeability ($K=1$) fracture channels (black lines) of width 0.1 embedded in a low-permeability background ($K=0.01$) in the domain $\Omega=[0,4]\times[0,1]$. A uniform grid with mesh size $h=1/64$ is used for the discretization.}
    \label{fig:fracture-setup}
\end{figure}

\begin{table}[]
    \centering
    \begin{tabular}{|c||c|c|c|}
    \hline
    Phase ($j$)  & Density ($\rho_j$) & Heat capacity ($c_j$) & Conductivity ($\lambda_j$)\\
    \hline 
     Incoming fluid ($i$)   &  1.0  & 1.0 & 1.0 \\
     resident fluid ($r$)   &  1.0  & 1.0 & 1.0 \\
     Solid ($s$)            &  2.0  & 2.0 & 2.0 \\
    \hline
    \end{tabular}
    \quad
    \begin{tabular}{|c|c|}
    \hline
    \multicolumn{2}{|c|}{Heat transfer coefficients} \\
    \hline
      $h_{ir}a_{ir}$   & 7.0 \\
      $h_{is}a_{is}$   & 1.0 \\
      $h_{rs}a_{rs}$   & 1.0 \\
    \hline
    \end{tabular}
    \quad
    \begin{tabular}{|c|l|}
    \hline
    \multicolumn{2}{|c|}{Other quantities} \\
    \hline
    Porosity              & $\phi = 0.5$ \\
    \hline
    \end{tabular}
    \caption{Example 2: parameters for phase properties, heat transfer, and other fixed quantities.}
    \label{tab:example2-parameters}
\end{table}

Figure \ref{fig:ex2-fluid-pressure} shows the steady-pressure profiles, and Figure \ref{fig:ex3-fluid-concentration} shows the concentration values at several time steps. We note that the concentration values remains to be bounded within $0.01 \leq z \leq 0.99$. 

Figure~\ref{fig:ex3-fluid-injected} shows injected-fluid temperature $T_i$ and its corresponding enthalpy $z T_i$. 
Although $T_i$ may attain relatively high values in regions where the injected-fluid concentration is low, the corresponding enthalpy remains small due to the concentration weighting.
The enthalpy carried by the injected fluid attains its maximum near the concentration fronts, where both temperature and concentration contribute significantly.

Similarly, Figure~\ref{fig:ex3-fluid-residing} shows the resident-fluid temperature $T_r$, along with its associated enthalpy. 
We observe that the resident-fluid enthalpy $(1-z)T_r$ is progressively extracted as heat is transferred to the injected fluid.

Finally, Figure~\ref{fig:ex3-fluid-Mix} displays the mixed-fluid temperature $T_f^{\text{mix}} = z T_i + (1-z)T_r$. 
The overall mixed-fluid temperature decreases over time due to continuous injection of cold fluid. For completeness, we show the solid temperature in Figure~\ref{fig:ex3-solid}, which remains relatively homogeneous because the heat transfer rate between the solid and fluid is slow. However, the solid temperature exhibits an overall decreasing trend with time.

\begin{figure}
    \centering
    \includegraphics[width=1.0\linewidth]{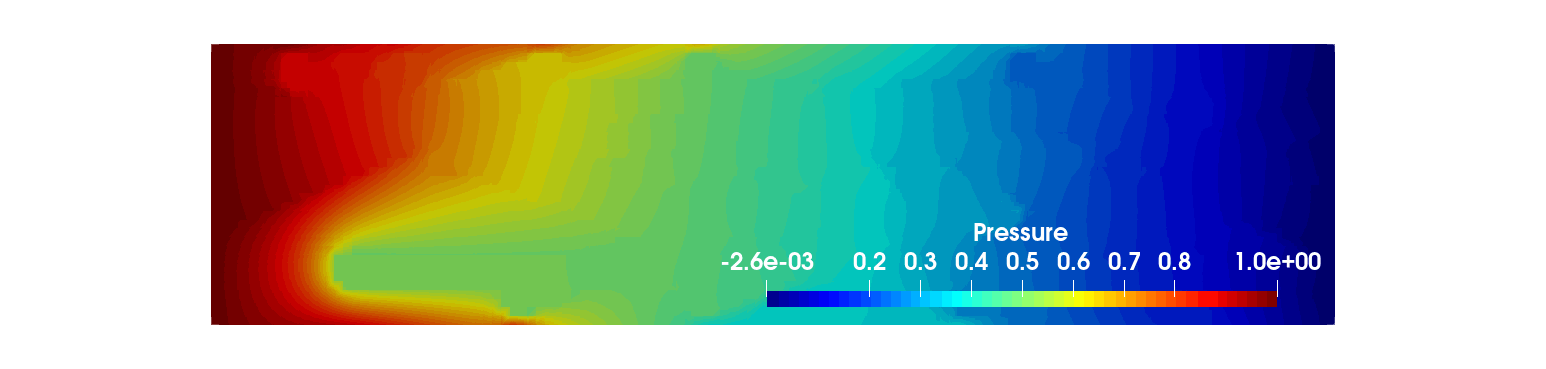}
    \caption{Example 3: pressure $p$.}
    \label{fig:ex2-fluid-pressure}
\end{figure}


\begin{figure}[H]
    \centering
    \begin{subfigure}[t]{1.0\textwidth}
    \centering\includegraphics[width=1.0\linewidth]{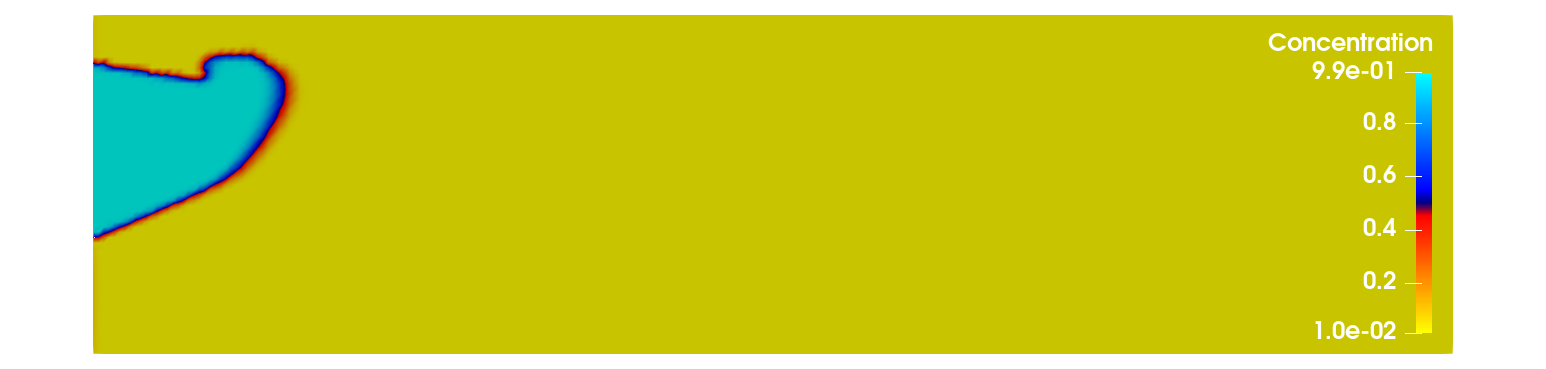}
    \caption{$t = 2.0$}
    \end{subfigure}
    \begin{subfigure}[t]{1.0\textwidth}
    \centering\includegraphics[width=1.0\linewidth]{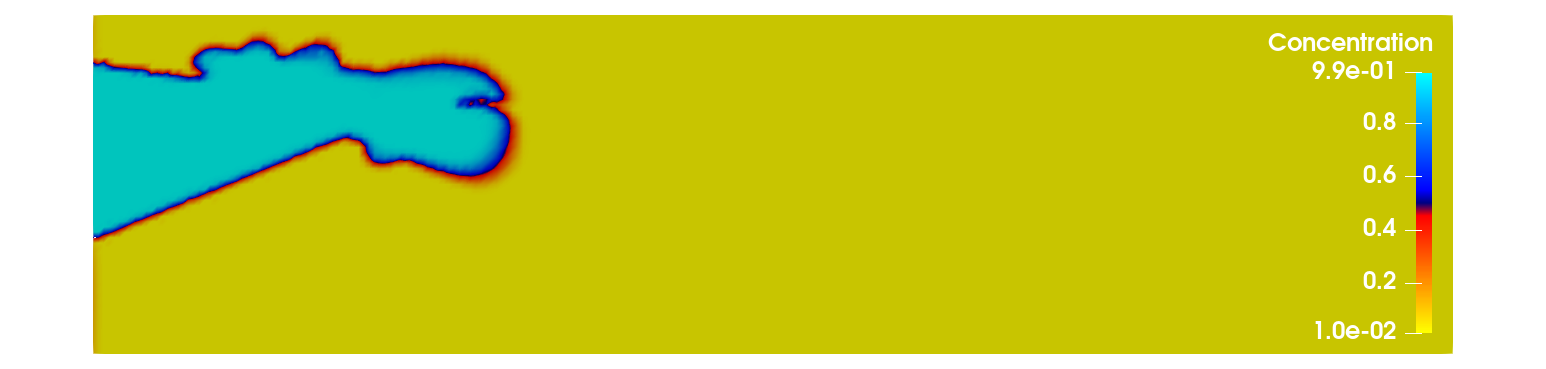}
    \caption{$t = 4.0$}
    \end{subfigure}
    \begin{subfigure}[t]{1.0\textwidth}
    \centering\includegraphics[width=1.0\linewidth]{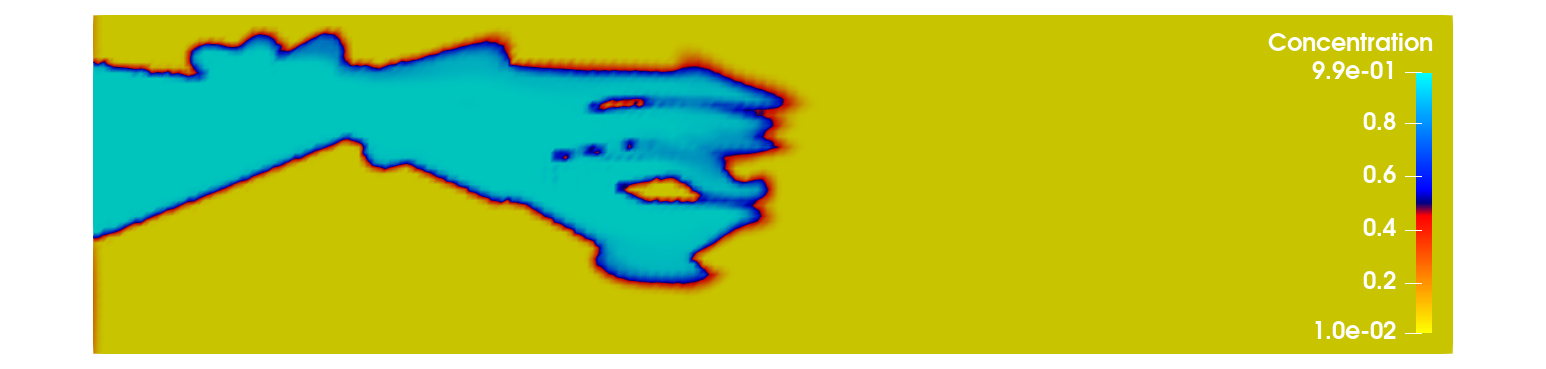}
    \caption{$t = 8.0$}
    \end{subfigure}
    \begin{subfigure}[t]{1.0\textwidth}
    \centering\includegraphics[width=1.0\linewidth]{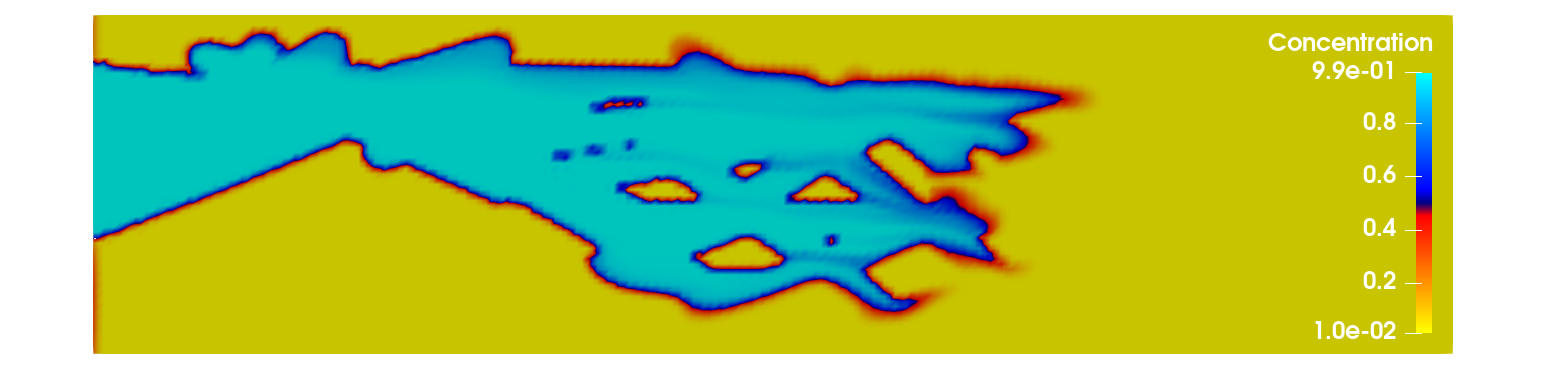}
    \caption{$t = 12.0$}
    \end{subfigure}
    \begin{subfigure}[t]{1.0\textwidth}
    \centering\includegraphics[width=1.0\linewidth]{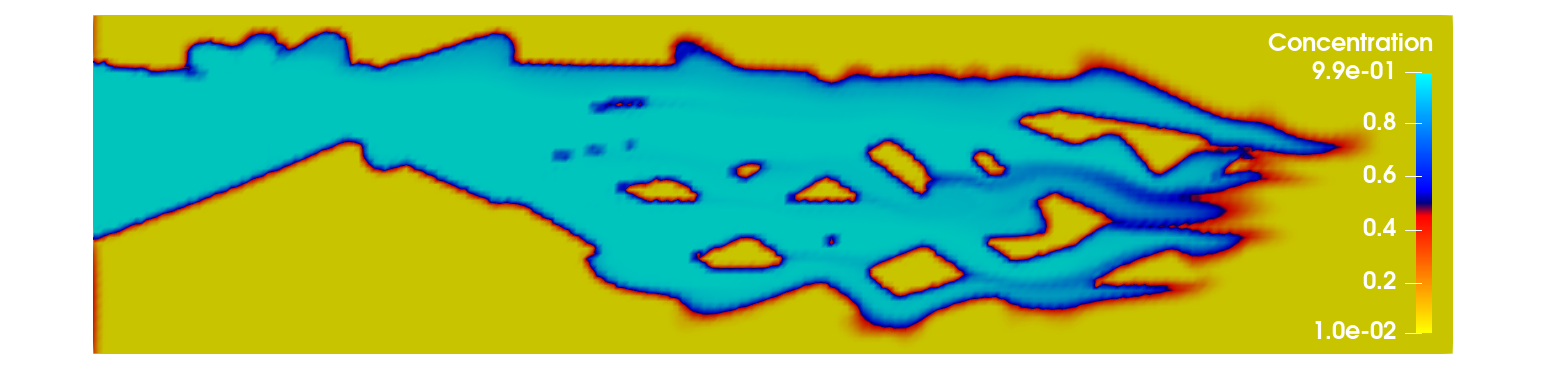}
    \caption{$t = 16.0$}
    \end{subfigure}
    \caption{Example 3: solution profiles for concentration $z$ at selected time step. }
    \label{fig:ex3-fluid-concentration}
\end{figure}



    

\begin{figure}
    \centering
    \begin{subfigure}[t]{1.2\textwidth}
    \centering
    \includegraphics[width=0.45\linewidth]{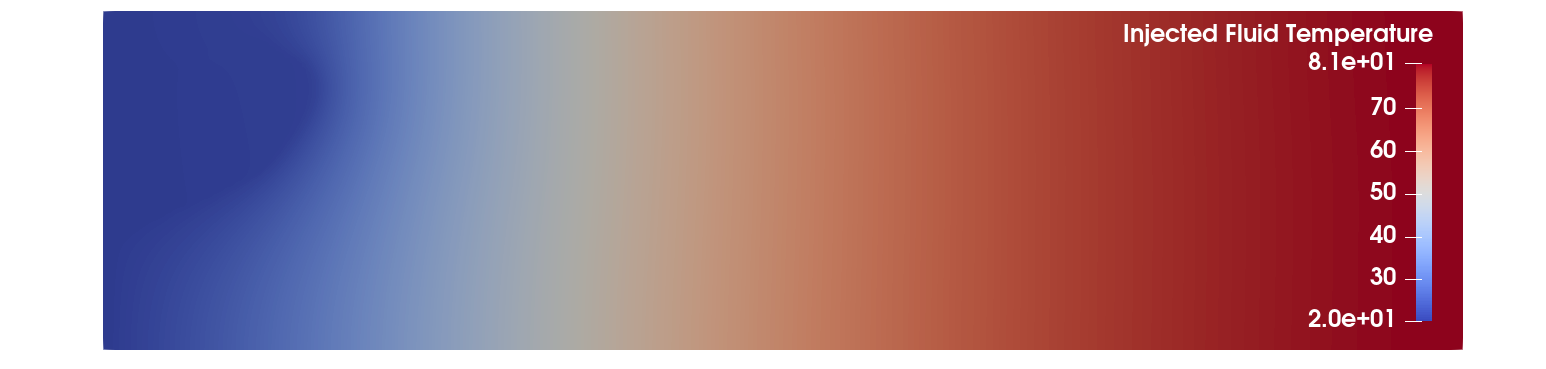}
    \includegraphics[width=0.45\linewidth]{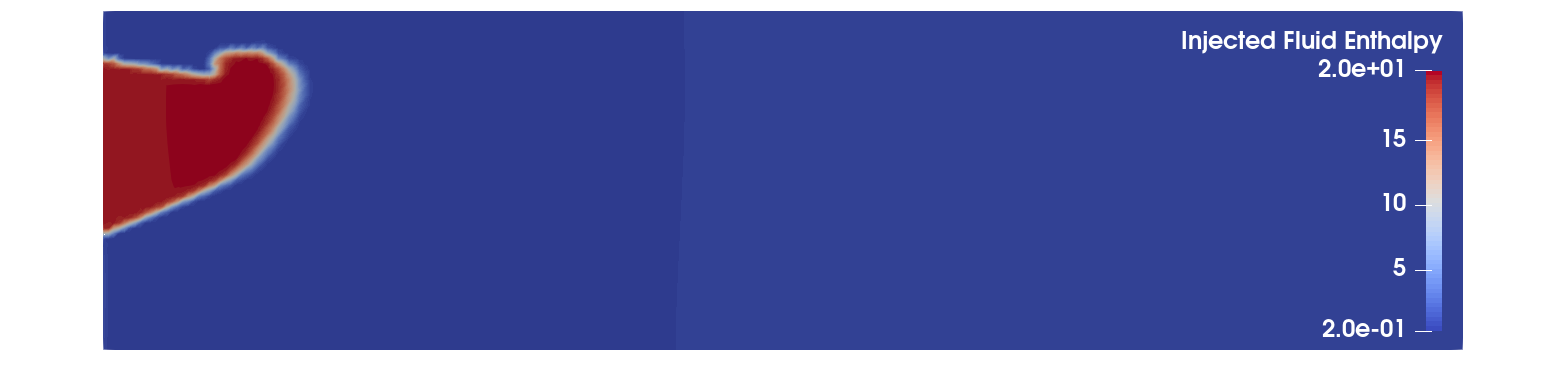}
    \caption{$t = 2.0$}
    \end{subfigure}

    \begin{subfigure}[t]{1.2\textwidth}
    \centering
    \includegraphics[width=0.45\linewidth]{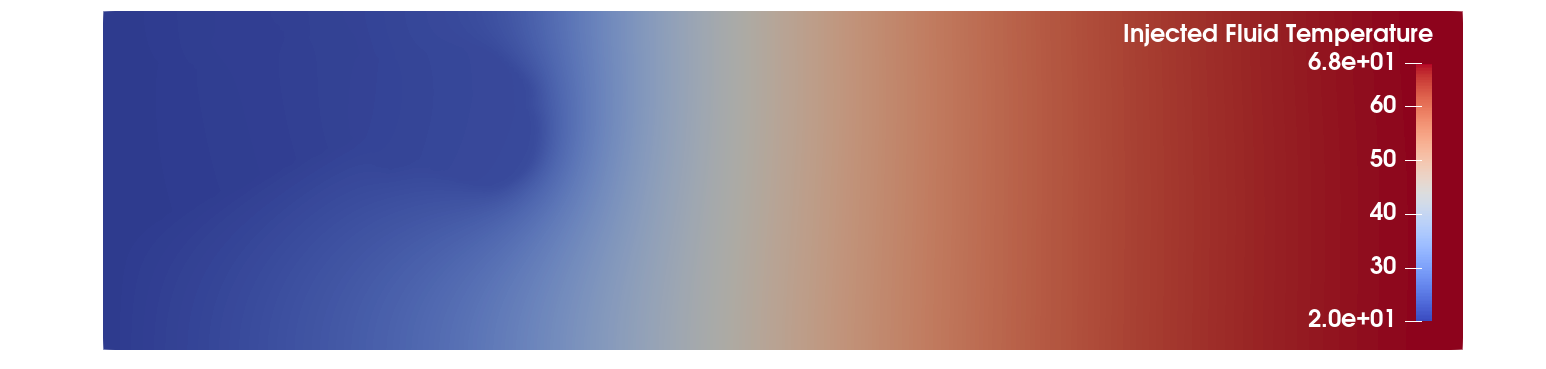}
    \includegraphics[width=0.45\linewidth]{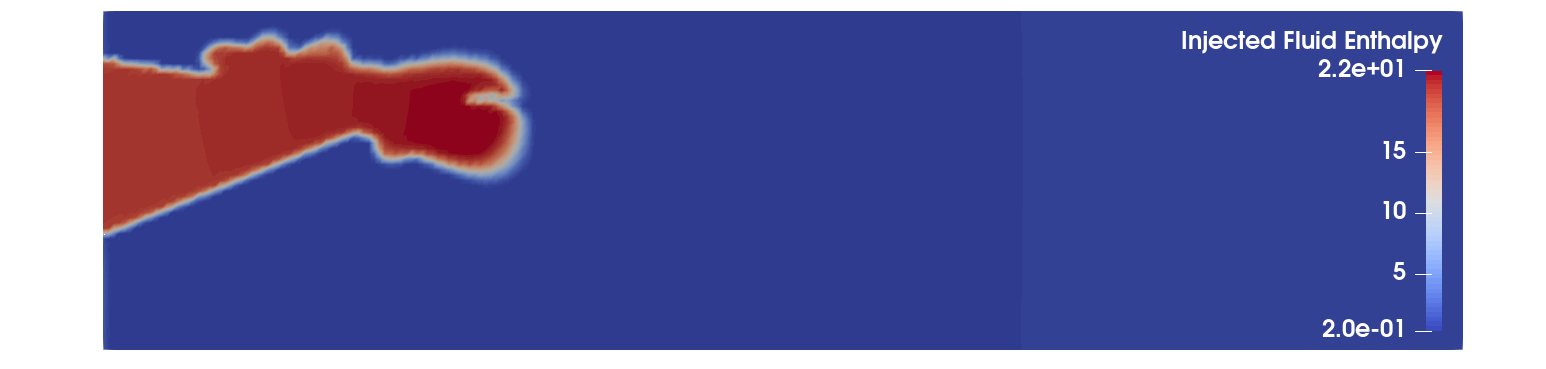}
    \caption{$t = 4.0$}
    \end{subfigure}

    \begin{subfigure}[t]{1.2\textwidth}
    \centering
    \includegraphics[width=0.45\linewidth]{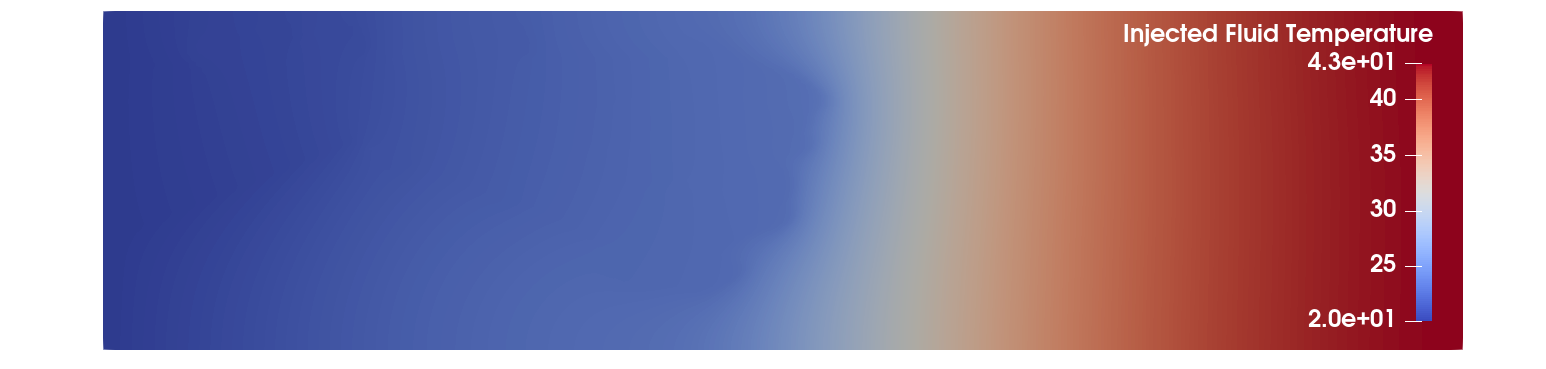}
    \includegraphics[width=0.45\linewidth]{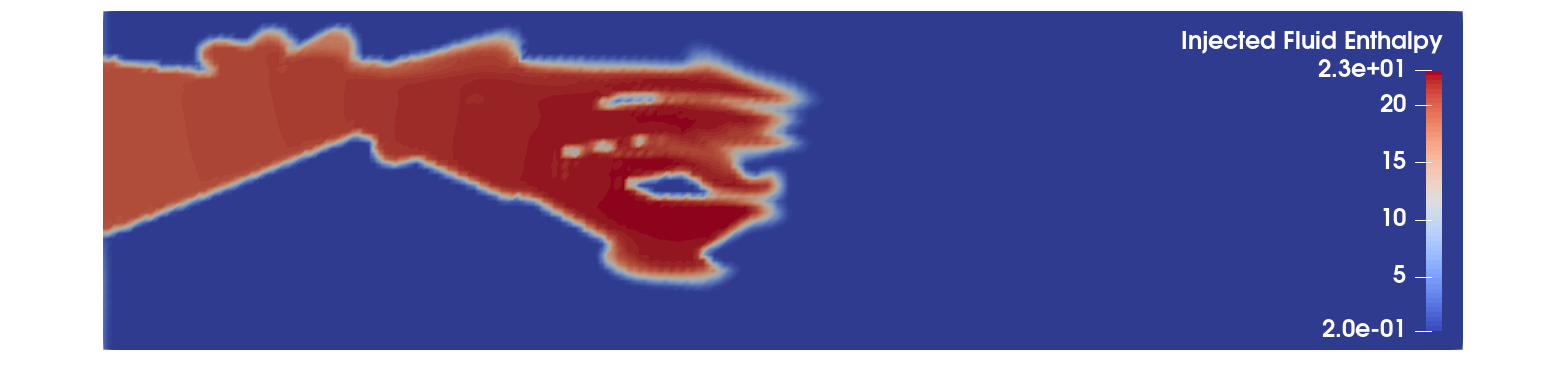}
    \caption{$t = 8.0$}
    \end{subfigure}

    \begin{subfigure}[t]{1.2\textwidth}
    \centering
    \includegraphics[width=0.45\linewidth]{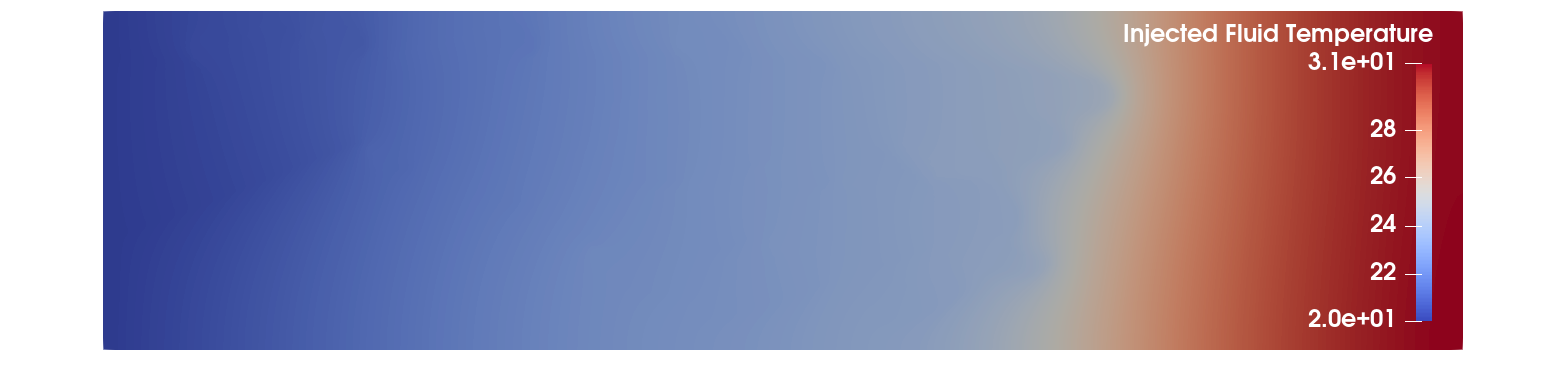}
    \includegraphics[width=0.45\linewidth]{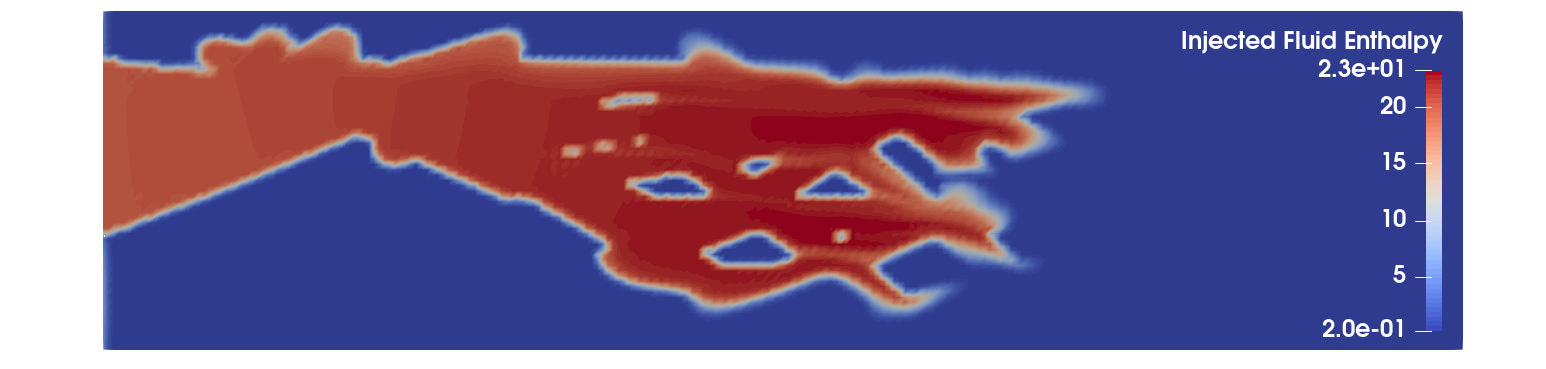}
    \caption{$t = 12.0$}
    \end{subfigure}

    \begin{subfigure}[t]{1.2\textwidth}
    \centering
    \includegraphics[width=0.45\linewidth]{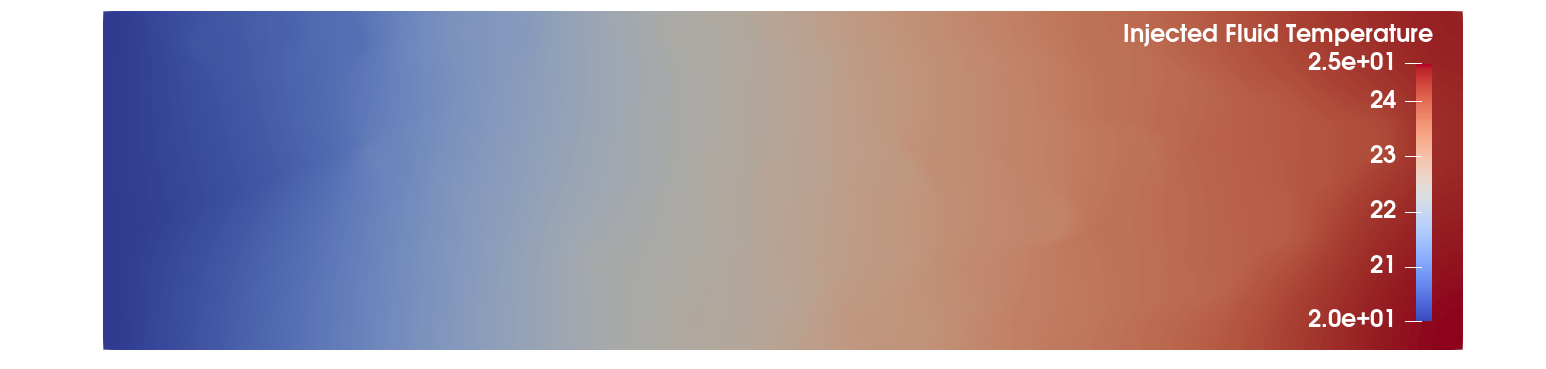}
    \includegraphics[width=0.45\linewidth]{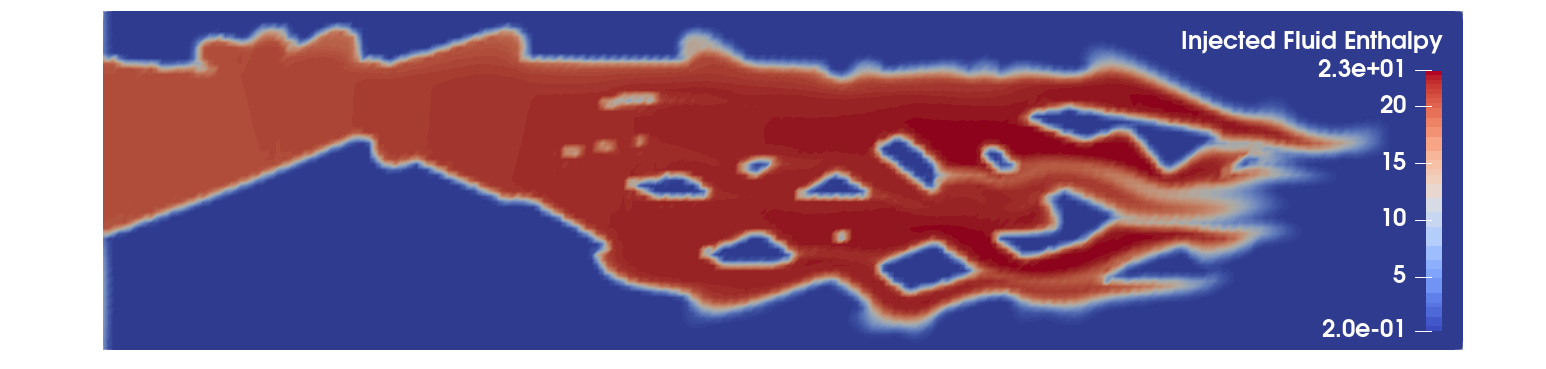}
    \caption{$t = 16.0$}
    \end{subfigure}
    
    \caption{Example 3: the solution profiles for injected fluid temperature $T_i$ (first column) and enthalpy $z\times T_i$ (second column) at selected time step. The colormap for enthalpy is discretized into 50 levels.}
    \label{fig:ex3-fluid-injected}
\end{figure}



    

\begin{figure}
    \centering
    \begin{subfigure}[t]{1.2\textwidth}
    \centering
    \includegraphics[width=0.45\linewidth]{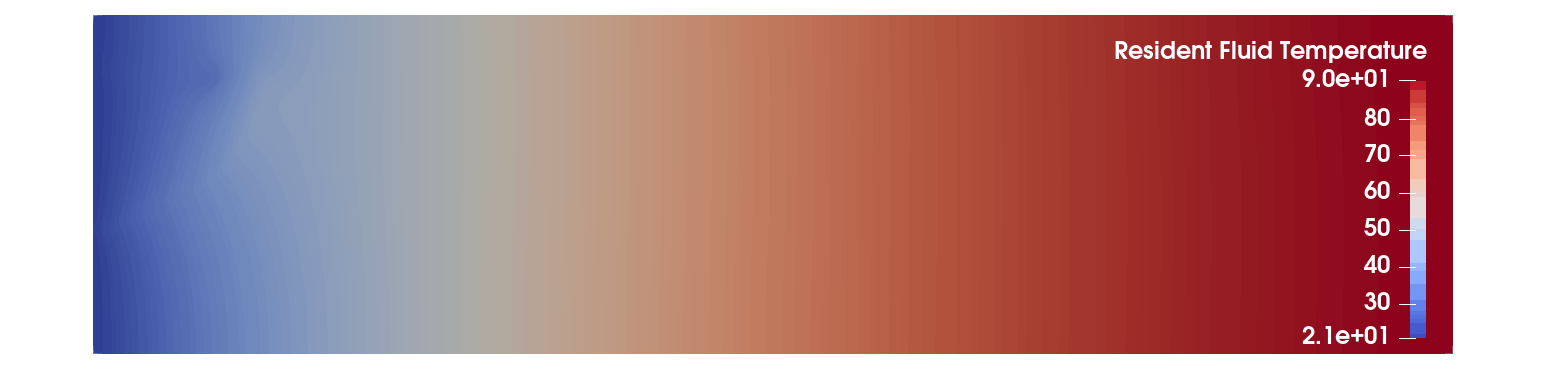}
    \includegraphics[width=0.45\linewidth]{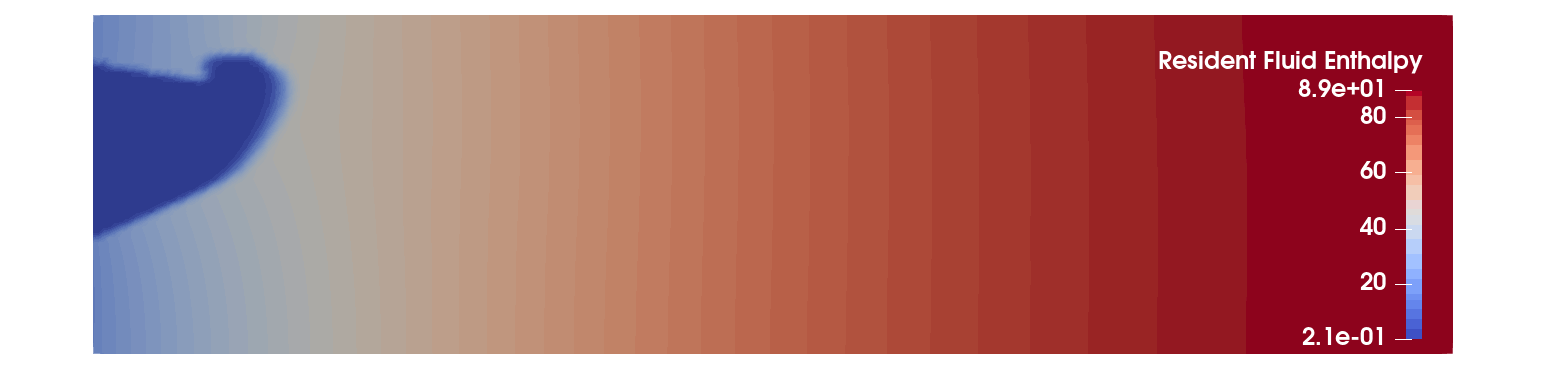}
    \caption{$t = 2.0$}
    \end{subfigure}

    \begin{subfigure}[t]{1.2\textwidth}
    \centering
    \includegraphics[width=0.45\linewidth]{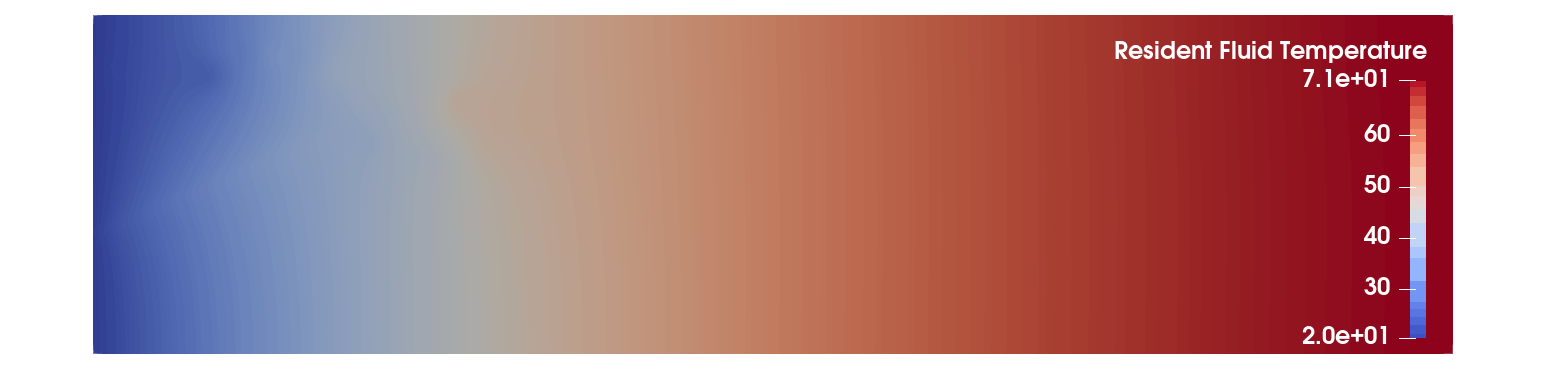}
    \includegraphics[width=0.45\linewidth]{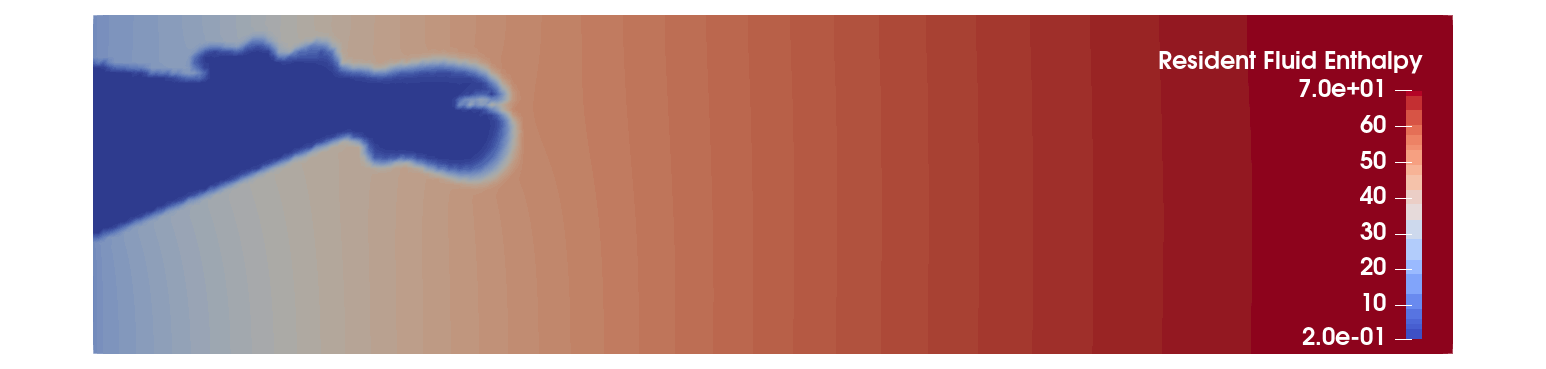}
    \caption{$t = 4.0$}
    \end{subfigure}

    \begin{subfigure}[t]{1.2\textwidth}
    \centering
    \includegraphics[width=0.45\linewidth]{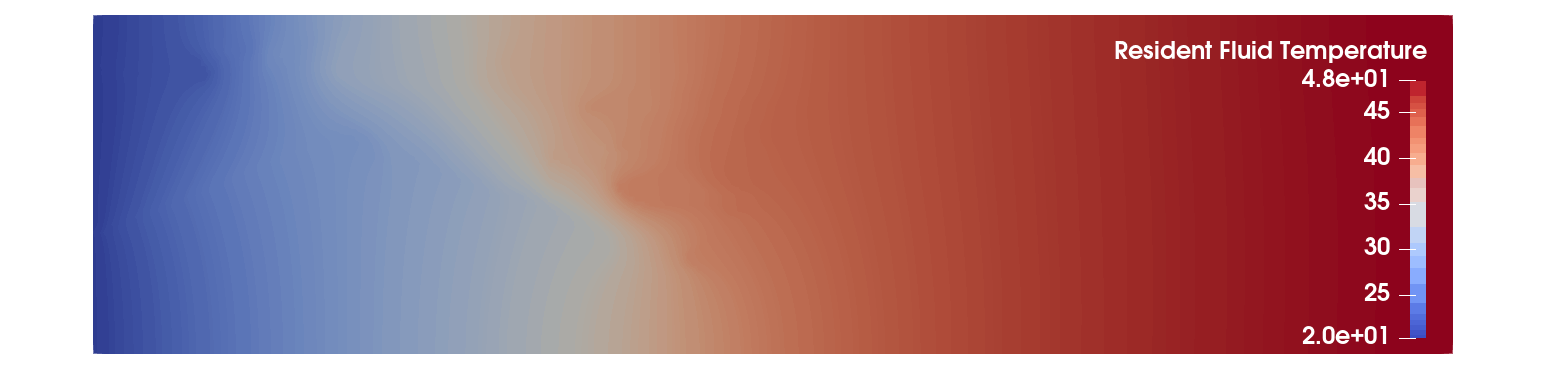}
    \includegraphics[width=0.45\linewidth]{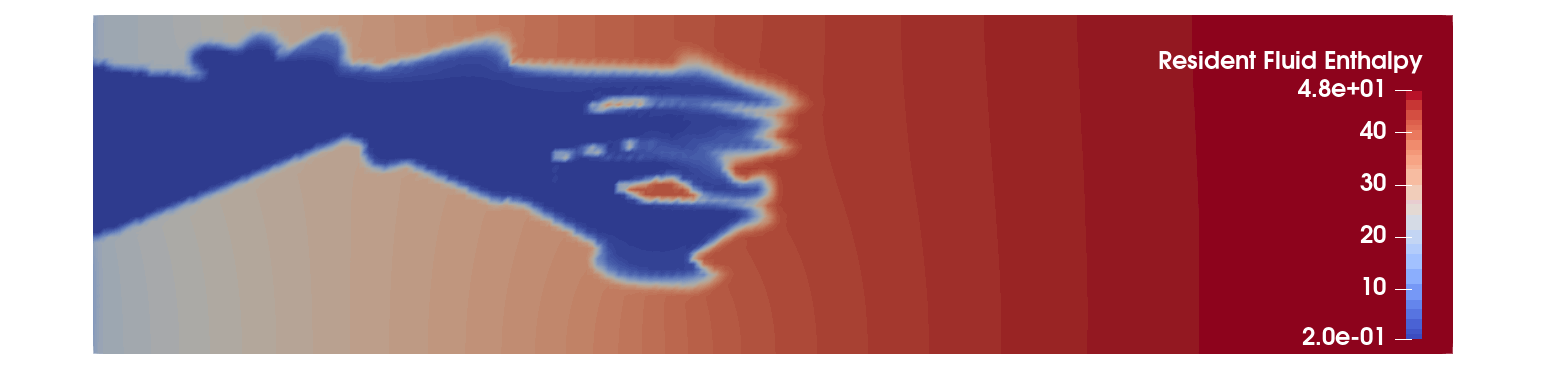}
    \caption{$t = 8.0$}
    \end{subfigure}

    \begin{subfigure}[t]{1.2\textwidth}
    \centering
    \includegraphics[width=0.45\linewidth]{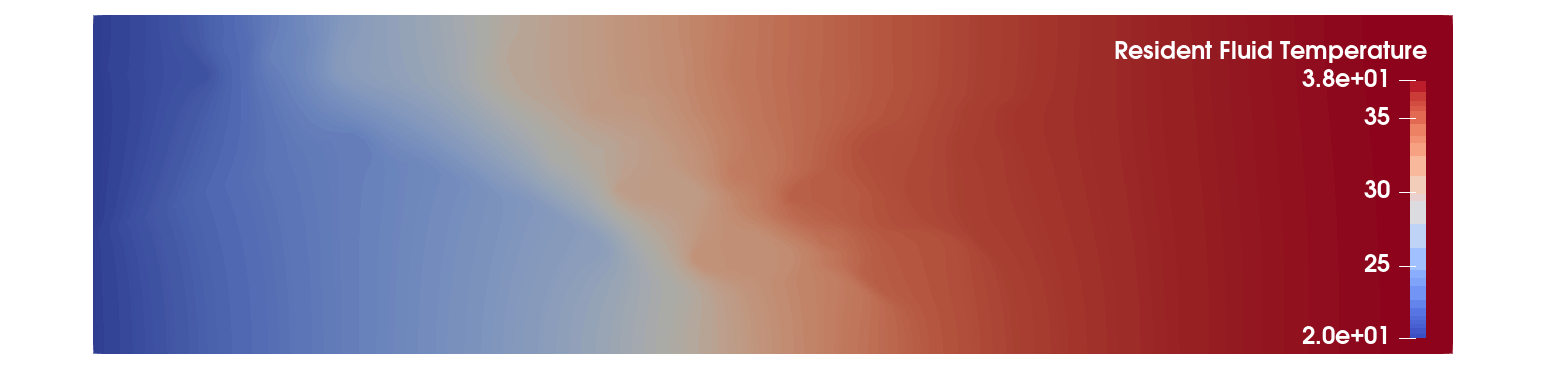}
    \includegraphics[width=0.45\linewidth]{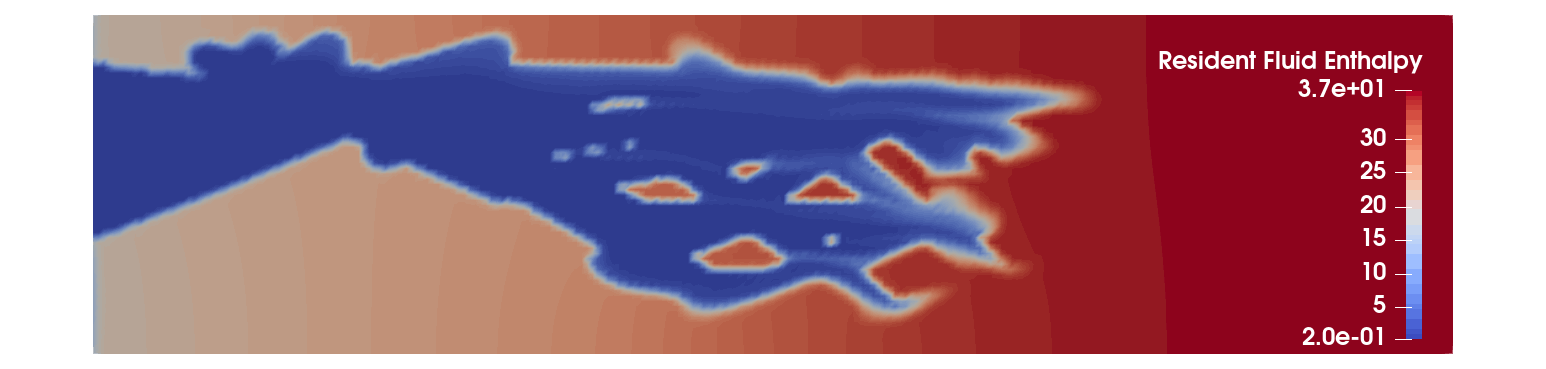}
    \caption{$t = 12.0$}
    \end{subfigure}

    \begin{subfigure}[t]{1.2\textwidth}
    \centering
    \includegraphics[width=0.45\linewidth]{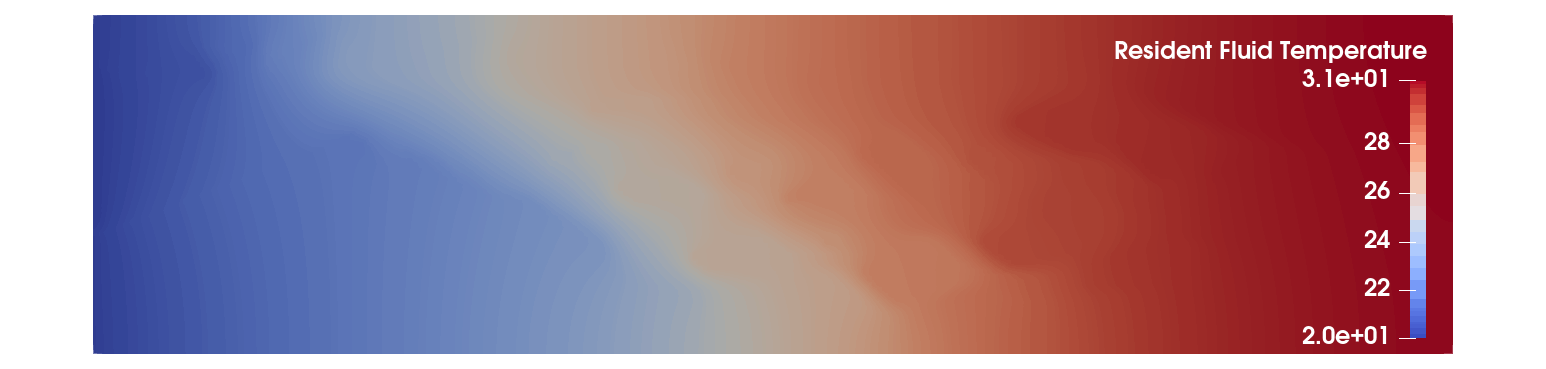}
    \includegraphics[width=0.45\linewidth]{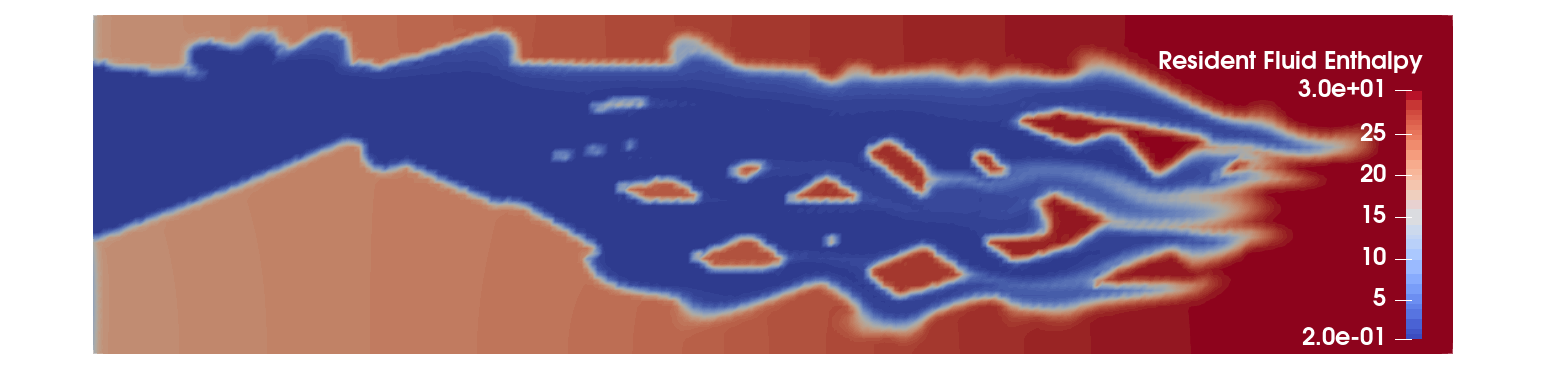}
    \caption{$t = 16.0$}
    \end{subfigure}
    
    \caption{Example 3: the solution profiles for injected fluid temperature $T_r$ (first column) and enthalpy $(1-z)\times T_r$ (second column) at selected time step. The colormap for enthalpy is discretized into 50 levels.}
    \label{fig:ex3-fluid-residing}
\end{figure}

\begin{figure}[H]
    \centering
    \begin{subfigure}[t]{1.0\textwidth}
    \centering\includegraphics[width=1.0\linewidth]{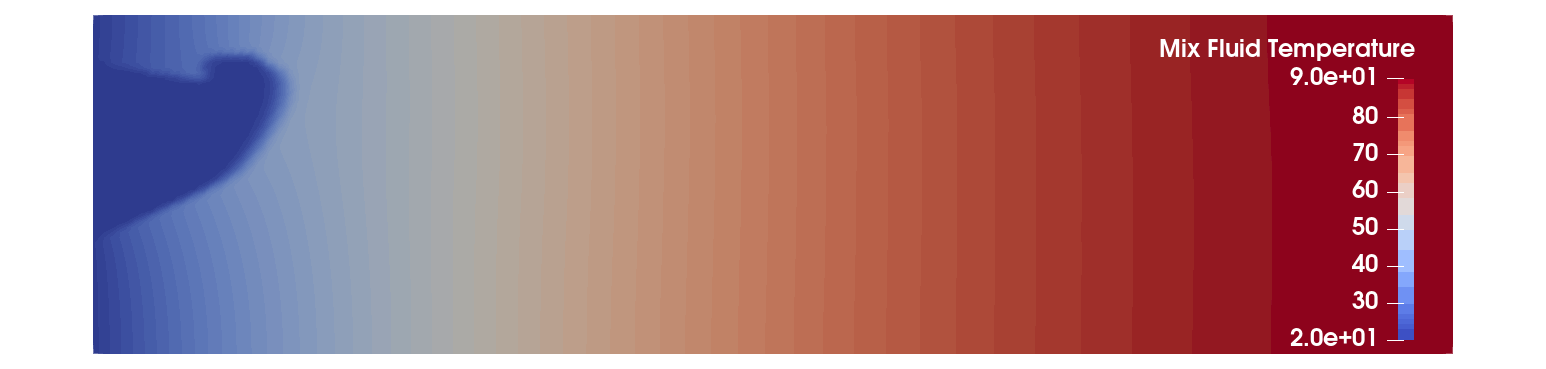}
    \caption{$t = 2.0$}
    \end{subfigure}
    \begin{subfigure}[t]{1.0\textwidth}
    \centering\includegraphics[width=1.0\linewidth]{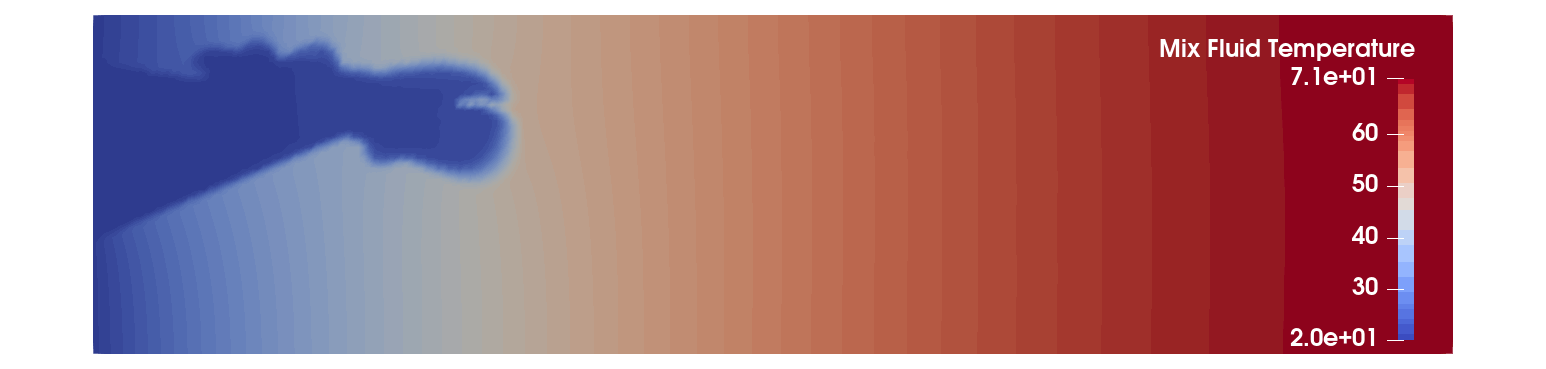}
    \caption{$t = 4.0$}
    \end{subfigure}
    \begin{subfigure}[t]{1.0\textwidth}
    \centering\includegraphics[width=1.0\linewidth]{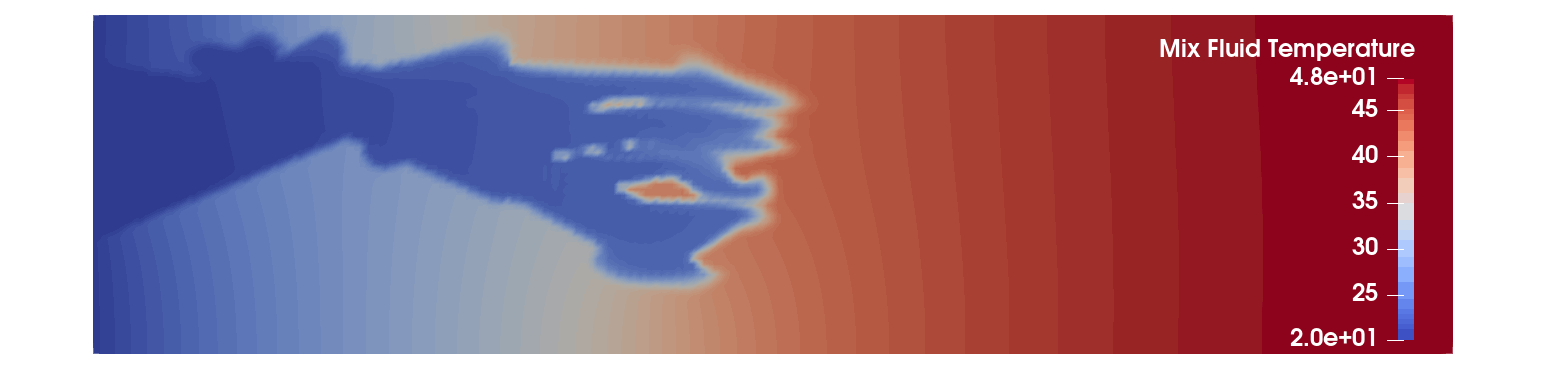}
    \caption{$t = 8.0$}
    \end{subfigure}
    \begin{subfigure}[t]{1.0\textwidth}
    \centering\includegraphics[width=1.0\linewidth]{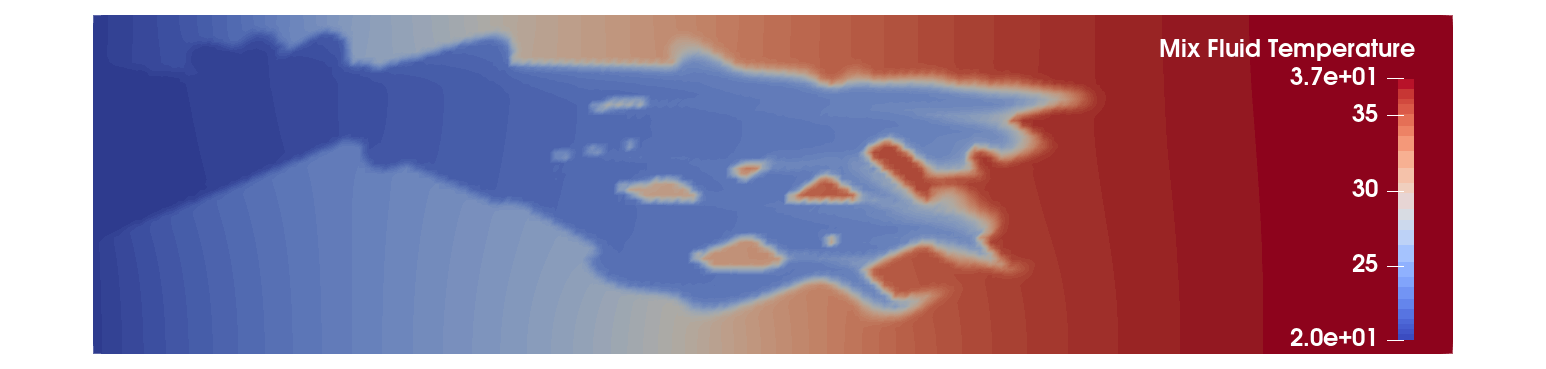}
    \caption{$t = 12.0$}
    \end{subfigure}
    \begin{subfigure}[t]{1.0\textwidth}
    \centering\includegraphics[width=1.0\linewidth]{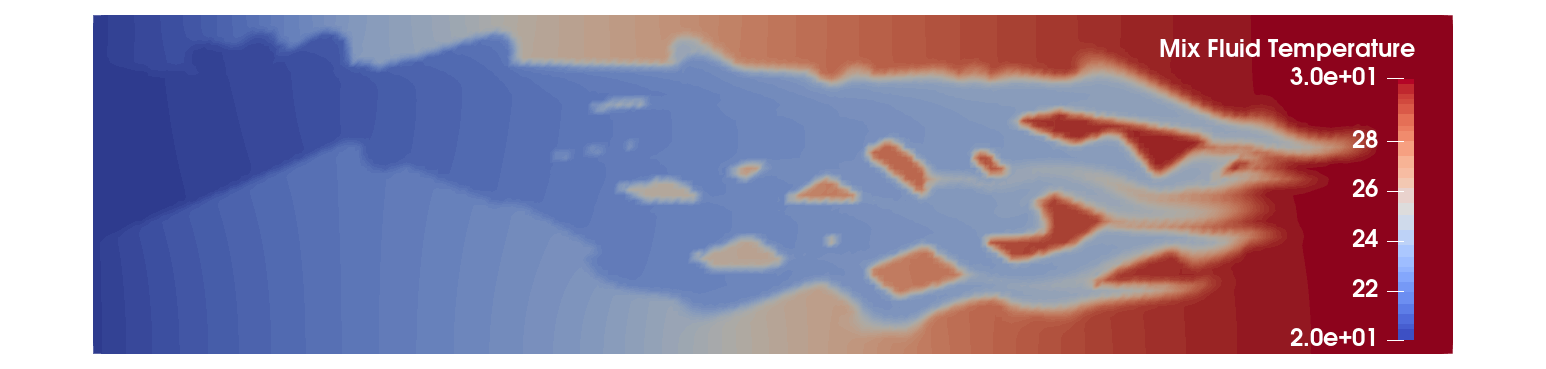}
    \caption{$t = 16.0$}
    \end{subfigure}
    \caption{Example 3: the profiles for mixed fluid temperature $T_f^\text{mix}:= zT_i + (1-z)T_r$ at selected time step. The colormap for mixed fluid temperature is discretized into 50 levels. }
    \label{fig:ex3-fluid-Mix}
\end{figure}

\begin{figure}
    \centering
    \begin{subfigure}[t]{1.0\textwidth}
    \centering
    \includegraphics[width=1.0\linewidth]{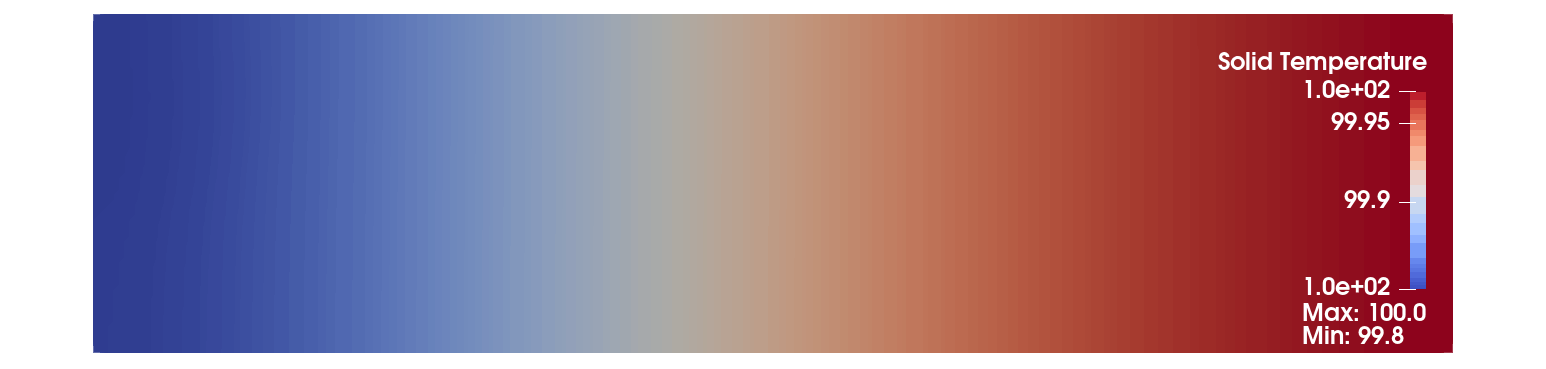}
    \caption{$t = 2.0$}
    \end{subfigure}
    \begin{subfigure}[t]{1.0\textwidth}
    \centering
    \includegraphics[width=1.0\linewidth]{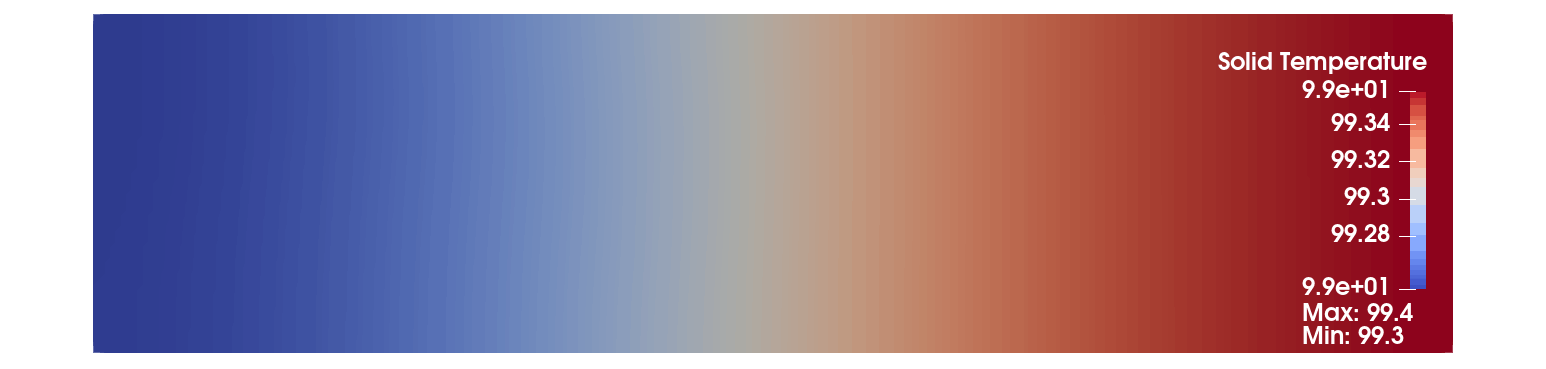}
    \caption{$t = 8.0$}
    \end{subfigure}
    \begin{subfigure}[t]{1.0\textwidth}
    \centering
    \includegraphics[width=1.0\linewidth]{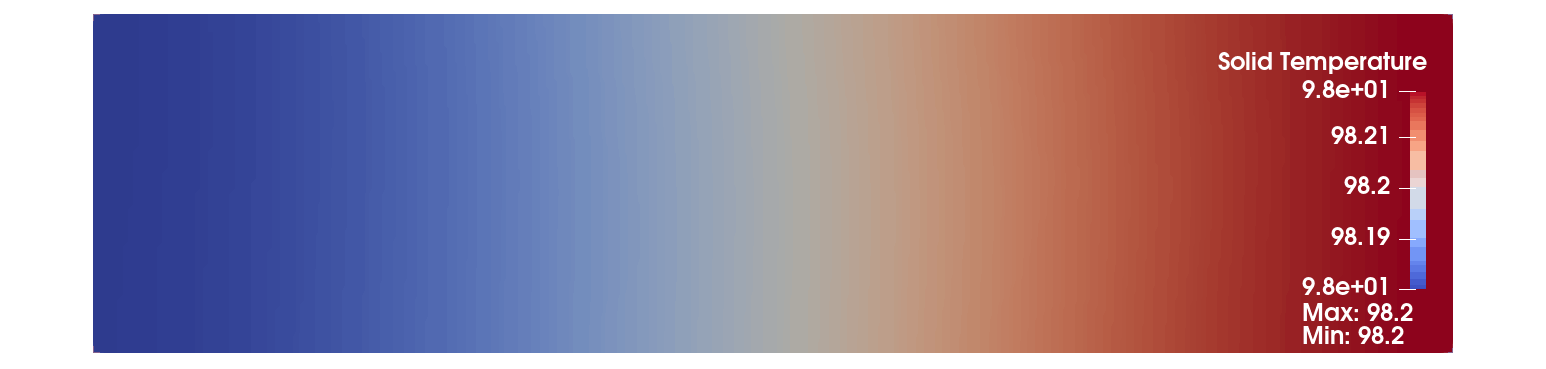}
    \caption{$t = 16.0$}
    \end{subfigure}
    \caption{Example 3: the solution profiles for solid temperature $T_s$ at selected time step.}
    \label{fig:ex3-solid}
\end{figure}

\paragraph{Conclusion}

In this work, we propose a three-way local thermal nonequilibrium (three-way LTNE) model that extends the classical LTNE framework by decomposing the fluid temperature into injected and resident components through a concentration field from a volumetric enthalpy perspective. 
The resulting formulation introduces separate temperature fields for injected fluid, resident fluid, and solid matrix, enabling explicit resolution of injected-fluid heating.

The model preserves thermodynamic consistency and reduces to classical LTNE and LTE formulations in appropriate limiting regimes.
An enriched Galerkin discretization combined with a sequential splitting strategy provides stable and locally conservative solutions for advection-dominated geothermal flow.

Numerical experiments demonstrate that the three-way LTNE model captures thermal features along concentration fronts that are not resolved by averaged-temperature approaches, highlighting its potential for improved prediction of thermal breakthrough in enhanced geothermal systems.

In the present study, the thermal diffusion coefficients are chosen sufficiently large to avoid strongly advection-dominated regimes. 
In convection-dominated settings, the BDF2 time discretization does not guarantee preservation of a discrete maximum principle, and additional stabilization or limiting strategies may therefore be required. 
Moreover, point sources are not considered in the pressure, concentration, or temperature equations. 
When such localized sources are present, further stabilization techniques may be necessary to prevent spurious oscillations in the numerical solution.

\bibliographystyle{elsarticle-num} 
\bibliography{reference}

\end{document}